\documentclass[12pt]{article}

\usepackage[cp866]{inputenc}      

\usepackage{amsfonts}
\usepackage{amssymb,amsmath}

\usepackage[T2A]{fontenc}    
\usepackage[english]{babel} 
\usepackage{epsfig}
\usepackage{graphicx}
\usepackage{amsmath}
\catcode`\@=11
\@addtoreset{equation}{section}\catcode`\@=12

\DeclareMathOperator{\arctanh}{arctanh}
\DeclareMathOperator{\spec}{spec}

\begin{document}

\newtheorem{lm}{Lemma}
\newtheorem{theorem}{Theorem}
\newtheorem{df}{Definition}
\newtheorem{prop}{Proposition}
\newtheorem{rem}{Remark}



\begin{center}
 {\large\bf Appearance of discrete Lorenz attractors in the transitions from saddle to saddle-focus}
\vspace{12pt}

 {\bf I.I. Ovsyannikov}
 \label{Author}
\vspace{6pt}

 {\small
 Constructor University, Germany; \\
 Lobachevsky State University of Nizhny Novgorod,
 Russia; \\

E-mail: ivan.i.ovsyannikov@gmail.com \\
iovsyannikov@constructor.university
  }
\end{center}

{\bf Abstract.} Strange attractors of the Lorenz type are among the most important
objects of this theory, because these chaotic attractors preserve their chaotic properties under
small perturbations. This feature is very important for applications, experiments, where certain
tolerance is unavoidably present, but it does not destroy chaotic attractors. The discrete
Lorenz attractor is a discrete-time analogue of the classical Lorenz attractor, but its dynamics
is richer - for example, it allows homoclinic tangencies, hence discrete Lorenz attractors are
wild attractors. Another important property of these attractors is that they can be born in 
bifurcations near triply degenerate fixed or periodic points. This makes possible analytic proof
of their existence in systems that allow the existence of such points.
Triply degenerate fixed points appear in global bifurcations - homoclinic and heteroclinic
tangencies. In order to get Lorenz-like attractors, the dynamics of the first return map along
the homoclinic or heteroclinic cycle should be effectively at least three-dimensional, i.e. there
should not exist lower-dimensional invariant manifolds. This can be achieved by adding some 
special conditions, global or local. Global degeneracies are related to the existence of non-
simple homoclinic tangencies or non-simple heteroclinic orbit in a cycle, these cases were
studied before. Local conditions either require the cycle to contain at least one saddle-focus,
or add certain relations on the multipliers of the fixed point such that the leading stable direction
of the saddle either disappears or alternates. All these cases were already studied before
except one, related to the transition from saddle to saddle-focus.
In the present paper this case is investigated, and the existence of a cascade of parameter
domains containing systems with discrete Lorenz attractors is proved. In particular, it includes
the Belyakov bifurcation, when the saddle becomes a saddle-focus through a collision of the
eigenvalues, and the 3DL bifurcation, when the dimension of the leading stable subset
alternates between 1 and 2. This paper completes the list of the simplest bifurcations of
homoclinic and heteroclinic tangencies by studying the last possible case.



\section*{Introduction}
%

Dynamical chaos is the research area that is being rapidly developed in the last decades. Chaotic phenomena are often represented by strange attractors. 
They can be of quite different types, and hyperbolic 
attractors and Lorenz attractors
are the mostly known and well investigated. Such attractors are remarkable by the fact that they preserve their
``strangeness'' under small smooth perturbations, and, besides, hyperbolic attractors are structurally stable\footnote{We note that many
well-known types of chaotic dynamics, such as H\'enon-like attractors, spiral attractors (R$\ddot{o}$ssler attractors, attractors observed in the
Chua circuits) etc, are, in fact, the so-called quasiattractors (in the terms introduced by Afraimovich and Shilnikov \cite{AS83})
because stable periodic orbits can appear in them under arbitrary small smooth perturbations.}.
Until the nearest time, only these attractors were known to compose the class of the so-called {\em genuine strange attractors}.
The situation was drastically changed after the paper \cite{TS98} by Turaev and Shilnikov  in which the theory of
the {\em wild hyperbolic attractors} was founded and an example of a wild spiral attractor was presented.

Wild hyperbolic attractors allow homoclinic tangencies, unlike the hyperbolic and Lorenz attractors, and therefore they belong to
 Newhouse domains \cite{NPT, GST93b}
but nevertheless stable periodic orbits as well as other stable invariant subsets are not born  under perturbations. 
Thus, wild hyperbolic attractors should also be referred to as genuine strange attractors\footnote{Moreover, by definition (see \cite{TS98}), such attractors are chain transitive and stable closed invariant sets. The chain transitivity means that any point of the attractor $\Lambda$ is admissible by $\varepsilon$-orbits from any other point of $\Lambda$; the stability means the existence of an open adsorbing  domain containing the attractor such that any orbit entering into the domain tends to $\Lambda$ (exponentially fast). Thus, this definition is, in fact, the definition of attractor by Ruelle and Conley.}.
An important subset of wild hyperbolic attractors is composed by discrete
Lorenz attractors which appear, in particular, in  Poincare maps for  periodically perturbed flows with the Lorenz attractors \cite{TS08}.
The definition of such an attractor is given in papers \cite{GGOT13, GGKS21}.
It is well-known that  Lorenz attractors do not allow homoclinic tangencies \cite{ABS77, ABS83}. However the latter can appear under non-autonomous periodic perturbations. Nevertheless,  bifurcation of these tangencies do not lead to the appearance of stable periodic orbits. 
The reason of it is that Lorenz-like attractors possess the pseudo-hyperbolic structure\footnote{This, very briefly, means that the differential $Df$of the corresponding map $f$ in the restriction onto the absorbing domain ${\cal D}$ of attractor $\Lambda$ admits an invariant splitting of form $E^{ss}_x\oplus E^{uc}_x$, for any point $x\in {\cal D}$, such that $Df$ is strongly contracting along directions $E^{ss}$ and expands volumes in transversal to $E^{ss}$ sections $E^{uc}$ (see \cite{TS98, TS08} for details)}.
However, for a long time, specific examples of systems of differential
equations with such attractors were not known. And quite recently, such an example
was proposed in \cite{GKT21}, this is a four-dimensional extension of the well-known Lorenz model. 

One of important peculiarities of the discrete Lorenz attractors is that such attractors can be born at local bifurcations
of periodic orbits having three or more
multipliers lying on the unit circle. This means that the corresponding attractors can be found in particular
models which have a sufficient number of parameters to provide the mentioned degeneracy. The following 3D H\'enon map
\begin{equation}
\bar x = y, \;\; \bar y = z, \;\; \bar z = M_1 + B x + M_2 y - z^2
\label{H3D}
\end{equation}
controlled by three independent parameters $M_1$, $M_2$ and $B$ is an example of such a model. 
In the papers \cite{GGOT13,GOST05,GMO06}
it was shown that the map (\ref{H3D}) possesses a discrete Lorenz-like attractor in some open parameter domain close to the point $(M_1=1/4, B=1, M_2 = 1)$ when the map has a fixed point with the triplet $(-1,-1,+1)$ of multipliers. Recently, discrete Lorenz attractors were also found in the orientation-reversing case $B < 0$, near point  $(M_1, M_2, B) = (7/4, -1, -1)$, see \cite{GGKS21}.

\begin{rem}
In paper \cite{GGS12} two universal scenarios of the development of chaos from a stable fixed point to
a strange attractor were represented for 3D maps. One of these scenarios leads to the appearance of a discrete Lorenz-like attractor.
This means that the discrete Lorenz attractors can be freely met in various dynamical models. In particular, such an attractor was found very recently in a
nonholonomic model of a rattleback \cite{GGK13}  and a model of two-component fluid convection \cite{EM18}.
\end{rem}

The fact that strange attractors of other nature than H\'enon-like ones (e.g. Lorenz-like attractors) can appear under bifurcations of homoclinic tangencies in multidimensional case (the dimension for maps is $\geq 4$) was announced first in the paper \cite{GST93c}. Concerning the
discrete Lorenz-like attractors, in the paper \cite{GMO06} it was shown that they appear in three parameter general
unfoldings in the case of three-dimensional diffeomorphisms having a homoclinic tangency to a fixed point
of the saddle-focus type with a Jacobian equal to one. Analogous results were obtained in papers \cite{GST09, GO10, GO13}
while studying three parameter general unfoldings of three-dimensional diffeomorphisms with a nontransversal
contracting-expanding heteroclinic cycle containing two fixed points which are saddle and saddle-focus \cite{GST09} or
are both saddle-foci \cite{GO10, GO13}.

Note that in all these cases it was assumed that at least one of the fixed points has the saddle-focus type.
Formally, together with the conditions on the Jacobians of the fixed points, it meant that the so-called {\em effective
dimension $d_e$ of the problem} (see \cite{T96}) equals to three i.e. is maximal for three-dimensional diffeomorphisms. When all fixed points are saddles, i.e. when all three multipliers of the fixed point are real, in the general case there exists a global two-dimensional invariant manifold, such that the transverse to it dynamics is trivial. Then, in order to keep the effective dimension at least three, additional degeneracies should be introduced that prevent from existence of such invariant manifolds. The full list of such conditions was identified in \cite{OvsPre}.

These degeneracies can be divided into global and local. Global degeneracies mean that the homoclinic and heteroclinic cycles contain non-simple orbits. Birth of discrete Lorenz attractors in such bifurcations were studied in paper \cite{GOT14} in the homoclinic case and 
paper \cite{O22} in the heteroclinic case. Local degeneracies occur when the eigenvalues at saddle fixed point satisfy resonance conditions. In paper \cite{GO17} the birth of discrete Lorenz attractors was studied when the stable eigenvalues are resonant, namely they satisfy condition $\lambda_1 = -\lambda_2$. Another possibility is to have multiple eigenvalues, such that in small perturbations the fixed point can become either a hyperbolic saddle or a saddle-focus. The continuous-time analog of this bifurcation was studied by L. Belyakov in paper \cite{Bel80}. So we call this bifurcation in discrete-time systems a discrete Belyakov bifurcation.

In the present paper we study cases of homoclinic and heteroclinic tangencies to fixed points having the saddle type, such that one of the fixed points undergoes the discrete Belyakov bifurcation.

The paper is organized as follows. In section~\ref{sec:def} there is the statement of the problem and the main result -- Theorem~\ref{thm:main} -- is formulated. In section~\ref{sec:first_ret} the first return map is constructed in both cases, and the non-degeneracy conditions are obtained. At the end of the section the rescaling lemma~\ref{th4} is formulated. The proof of the lemma is in section~\ref{sec:th3proof}.


\section{Statement of the problem and main results}\label{sec:def}

We consider two cases of $C^r$-diffeomorphisms, $r \geq 3$, satisfying the following conditions:

{\bf Case I: A homoclinic Belyakov bifurcation}  

{\bf I.A)} $f_0$ is three-dimensional;

{\bf I.B)} $f_0$ has a fixed point $O$ with multipliers $(\lambda_1, \lambda_2, \gamma)$ where
    $ \lambda_1 = \lambda_2 = \lambda$ and  $|\lambda| < 1 < |\gamma|$;

{\bf I.C)} The Jacobian $J(O) \equiv |\lambda^2 \gamma|$ at point $O$ is equal to one;

{\bf I.D)} Unstable manifold $W^u(O)$ and stable manifold $W^s(O)$ have a quadratic tangency at the points of a homoclinic orbit $\Gamma_0$.
\\

{\bf Case II: A heteroclinic cycle with Belyakov bifurcation}

{\bf II.A)} $f_0$ is three-dimensional;

{\bf II.B)} $f_0$ has two fixed points $O_1$ and $O_2$ with multipliers $(\lambda_{1}, \lambda_2, \gamma_1)$ and $(\nu_1, \nu_2, \gamma_2)$ respectively, such that $ \lambda_1 = \lambda_2 = \lambda$,  $|\lambda| < 1 < |\gamma_1|$,
$|\nu_2| < |\nu_1| < 1 < |\gamma_2|$;

{\bf II.C} One of the Jacobians $J_1 = |\lambda_1 \lambda_2 \gamma_1|$ and 
$J_2 = |\nu_1 \nu_2 \gamma_2|$ is greater than $1$, and another one is smaller than $1$;

{\bf II.D.a} Unstable manifold $W^u(O_1)$ and stable manifold $W^s(O_2)$ intersect transversely at the points of heteroclinic orbit $\Gamma_{12}$, 
unstable manifold $W^u(O_2)$ and stable manifold $W^s(O_1)$ have quadratic tangencies at the points of heteroclinic orbit $\Gamma_{21}$;

{\bf II.D.b} Unstable manifold $W^u(O_1)$ and stable manifold $W^s(O_2)$ have quadratic tangencies at the points of heteroclinic orbit $\Gamma_{12}$, 
unstable manifold $W^u(O_2)$ and stable manifold $W^s(O_1)$ intersect transversely at the points of heteroclinic orbit $\Gamma_{21}$. \\






%
%

We will study bifurcations of single-round periodic orbits in generic parametric unfoldings $f_\mu$ such that
$\left. f_\mu \right|_{\mu = 0} = f_0$. The main attention will be paid to the birth of discrete Lorenz attractors in local bifurcations. 
First, we have to identify the necessary set of parameters required to unfold the original system. 

\begin{figure}
\centerline{\epsfig{file=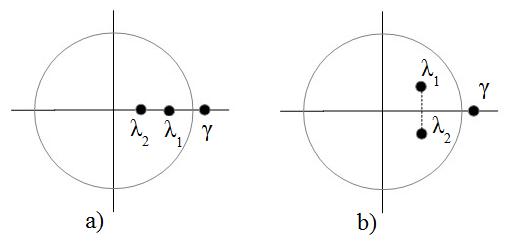,
width=8cm,
}} \caption{{\footnotesize The multipliers in the fixed point $O$ in Case I, or $O_2$ in Case II a) for $\mu_2 > 0$; b) for $\mu_2 < 0$.}} \label{3D}
\end{figure}

In Case I we take, as usual, the splitting distance of manifolds $W^u(O)$ and $W^s(O)$ with respect of some of the points
of $\Gamma_0$ as the first parameter $\mu_1$. At the bifurcation moment we have $\lambda_1 = \lambda_2$ (condition {\bf I.B}), and we
introduce parameter $\mu_2$ to unfold this condition, such that for $\mu_2 = 0$ the multipliers are equal $\lambda_1(\mu) = \lambda_2(\mu)$, for $\mu_2 > 0$ both eigenvalues are real
and for $\mu_2 < 0$ they are a complex-conjugate pair, see fig.~\ref{3D}.

The third parameter will control the Jacobian at $O$ in the way that
\begin{equation} \label{mu3.I}
\mu_3 = 1 - |\lambda_1(\mu) \lambda_2(\mu) \gamma(\mu)|.
\end{equation}


In Case II value $\mu_1$ will be also a splitting parameter of a quadratic tangency, $\mu_2$ will also control the multipliers of the fixed point $O_1$ such that they are real for $\mu_2 > 0$ and complex for $\mu_2 < 0$. Also note that due to condition {\bf II.B}, near one saddle point the phase volumes are contracted and near another one the volumes are expanded. Thus, we are dealing with the contracting-expanding (or mixed) case, which requires one more 
parameter $\mu_3$, that will control the values of the Jacobians 
$J_1$ and $J_2$. Similar to \cite{GST09, GO10, GO13, O22}, the following value will play this role:
\begin{equation}
\label{eq:e3}
\mu_3 = S(f_\mu) - S(f_0), 
\end{equation}
where $S(f)$ is a functional defined as $\displaystyle S(f) = - \frac{\ln J_1}{\ln J_2}$.

In any neighbourhood 
of initial diffeomorphism $f_0$ there exist diffeomorphisms $f_\mu$ having saddle fixed points
with non-trivial leading and non-leading stable invariant subspaces. 
This means that homoclinic or heteroclinic orbits to such fixed points can potentially have 
additional degeneracies. In order to prevent from other bifurcations, we have to add generality conditions. 
Namely, we require the homoclinic or heteroclinic orbits to be {\em simple} in some neighbourhood of $f_0$. 

The main result of the present paper is the following theorem.

\begin{theorem}[Main Theorem] \label{thm:main}
Let $f_\mu$ be the three-parametric family under consideration, that is, $f_0$ satisfies  A)--D), and $f_\mu$ is a general unfolding under conditions B), C) and D). Then, 

1) In the homoclinic Case I
in any neighbourhood of the origin $\mu = 0$ in the parameter space there exist infinitely many domains
$\delta_k$, where $\delta_k\to (0,0,0)$ as $k\to\infty$, such that
the diffeomorphism $f_\mu$  has a discrete Lorenz attractor at $\mu\in\delta_k$. 

2) In the heteroclinic Case II in any neighbourhood of the origin $\mu = 0$ in the parameter space there exist infinitely many domains
$\delta_{kj}$, where $\delta_{kj}\to (0,0,0)$ as $k, j \to \infty$, such that
the diffeomorphism $f_\mu$  has a discrete Lorenz attractor at $\mu\in\delta_{kj}$.
\end{theorem}

The proof of the Theorem is based on the Rescaling Lemma~\ref{th4}, and is given after the lemma's formulation.



\section{Construction of the first return map}\label{sec:first_ret}
The first return maps near homoclinic or heteroclinic cycles is a composition of local and global maps. The local maps are defined in some small neighbourhoods of fixed points, and the global maps connect these neighbourhoods along the global pieces of homoclinic or heteroclinic orbits. 

We consider a sufficiently small neighbourhood $U$ of the homoclinic of the heteroclinic cycle. It is composed 
as a union of a small neighborhood $U_{0}$ of fixed point $O$ with a finite number of small neighborhoods $U_i$ of those
points of $\Gamma$ which do not belong to $U_{0}$ in Case I. In the heteroclinic Case II, it is composed of small neighbourhoods $U_1$ and $U_2$ of fixed points $O_1$ and $O_2$ and a finite number of small neighborhoods $U_i$ of those
points of $\Gamma_{12}$ and $\Gamma_{21}$ which do not belong to $U_{1,2}$.
Each single-round periodic orbit lying entirely in $U$
has exactly one intersection point with each of $U_i$.

The local maps $T_0, T_{10}, T_{20}$ near respective fixed points $O, O_1, O_2$ is defined as a restriction of $f_\mu$ onto $U_0$, $U_1$ and $U_2$. 

We select two homoclinic points $M^+ \in U_0 \cap \Gamma$ and $M^- \in U_0 \cap \Gamma$, such that $M^+ \in W^s_{loc}(O)$ and $M^- \in W^u_{loc}(O)$ and their respective small neighbourhoods $\Pi^+ \subset U_0$ and $\Pi^- \subset U_0$. It is clear that there exists a natural number $q$ such that $M^+ = f_0 M^-$. We define the global map as $T_1(\mu)= f_\mu^q: \Pi^- \to \Pi^+ $.

In the heteroclinic case, we select in $U_1$ heteroclinic points $M_1^+ \in W^s_{loc}(O_1)$ and $M_1^- \in W^u_{loc}(O_1)$ and in $U_2$ heteroclinic points $M_2^+ \in W^s_{loc}(O_2)$ and $M_2^- \in W^u_{loc}(O_2)$, such that $M_1^-, \; M_2^+ \in \Gamma_{12}$ and $M_2^-, \; M_1^+ \in \Gamma_{21}$. There exist numbers $q_1$ and $q_2$ such that $M_2^+ = f_0^{q_1} M_1^-$ and $M_1^+ = f_0^{q_2} M_2^-$. We define global maps along heteroclinic orbits $\Gamma_{12}$ and $\Gamma_{21}$ as $T_{12} = f_\mu^{q_1} : \Pi_1^- \to \Pi_2^+$ and
$T_{21} = f_\mu^{q_2} : \Pi_2^- \to \Pi_1^+$


Begin iterating $\Pi^+$ under the action of $T_{0}$. Starting from some number $k_0$ images
$T_{0}^{k} \Pi^+$ will have a nonempty intersection with $\Pi^-$. The first return map $T_{k} \equiv T_{1} T_{0}^k$ is thus defined for every $k \ge k_0$ on an infinite set of stripes
$\sigma_{0}^k \subset \Pi^+$ accumulating to $W^s(O)$ as $k \to \infty$. Their
images under the iterations of the local map are stripes $f^k_\mu \sigma_{0}^k = \sigma_{1}^k \subset \Pi^-$
which accumulate to $W^u(O)$. 
Single-round periodic orbits lie in intersections $\sigma_{0}^k \cap T_1 \sigma_{1}^k$.

Analogously to the homoclinic case, in the heteroclinic case there exists a number $k_0$ such that images $T_{10}^k \Pi_1^+$ have a non-empty intersection with $\Pi_1^-$ for $k > k_0$ and respectively $j_0$ such that $T_{20}^j \Pi_2^+$ have a non-empty intersection with $\Pi_2^-$ for $j > j_0$. The first return maps $T_{kj} = T_{21}T_{20}^jT_{12}T_{10}^k$ is therefore defined on subsets $V_{kj} \subset \Pi_1^+$ that accumulate to $W^s(O_1)$ as $k, j \to \infty$.
In order to write explicit formulas for the first return map we need to write both local and global maps
in the most suitable form.

%

\subsection{Properties of the local maps}\label{sect:local}

 Assume that 
a fixed point possesses $n$ leading stable eigenvalues $\lambda_1, \ldots, \lambda_{n}$ and $m$ non-leading stable eigenvalues $\lambda_{n + 1}, \ldots, \lambda_{n + m}$ and an unstable eigenvalue $\gamma$. This means that $|\lambda_1| = \ldots = |\lambda_n| = \lambda < 1$ and $|\lambda_i| < \lambda$ for $i > n$.  
The main normal form has then the following representation in some $C^r$-smooth local
coordinates (see \cite{GS90, GST07, GS92, book}):
\begin{equation}
\begin{array}{l}
\bar x_1 \; = \; A_1(\mu) x_1 + \tilde H_2(y, \mu)x_2 + O(\|x\|^2\|y\|)  \\
\bar x_2 \; = \; A_2(\mu) x_2 + \tilde R_2(x, \mu) + \tilde H_4(y, \mu)x_2 + O(\|x\|^2\|y\|) \\
\bar y \; = \; \gamma(\mu) y + O(\|x\| \|y\|^2),\\
\end{array}
\label{t0norm}
\end{equation}
%
%
where $\spec A_1(0) = \{ \lambda_1, \ldots, \lambda_{n}\}$, $\spec A_2(0) = \{ \lambda_{n + 1}, \ldots, \lambda_{n + m}\}$,
$\tilde H_{2, 4}(0, \mu) = 0$, $\tilde R_{2}(x, \mu) = O(\|x\|^2)$.
In coordinates (\ref{t0norm}) the invariant manifolds of saddle fixed point $O$ are 
locally straightened: stable $W^s_{loc}(O): \{ y = 0\}$,
unstable $W^u_{loc}(O): \{ x_1 = 0, \; x_2 = 0\}$ and strong stable $W^{ss}_{loc}(O):  \{ x_1 = 0, \; y = 0\}$.

If the fixed points is a saddle ($O_2$), then matrices $A_1(\mu)$ and $A_2(\mu)$ are one-dimensional, $\spec A_1 = \{ \nu_1 \}$ and $\spec A_2 = \{ \nu_2 \}$. Multiple iterations of the local map $T_{20}^j :\; U_{2} \to U_{2}$ near such a point can be represented in the form close to linear (see, for example, \cite{GS92,book}):
\begin{equation}\label{eq:T20j}
\begin{array}{l}
  u_{1j} = \nu_1^j(\mu) u_{10} +
            \hat\nu^{j} \eta_{1j}(u_0, v_k, \mu)\\
  u_{2j} = \hat\nu^{j} \eta_{2j}(u_0, v_k, \mu)\\
  v_0  = \gamma_2^{-j} v_j + \hat \gamma_2^{-j} \eta_{3j}(u_0, v_j, \mu).
\end{array}
\end{equation}
Here $0 < \hat \nu < \nu_1$, $\hat\gamma_2 > \gamma_2$, 
 functions $\eta_{j}$ are uniformly bounded in $j$ together with their derivatives up to order $(r - 2)$,
 and their higher-order derivatives tend to zero.

\begin{figure}
\centerline{\epsfig{file=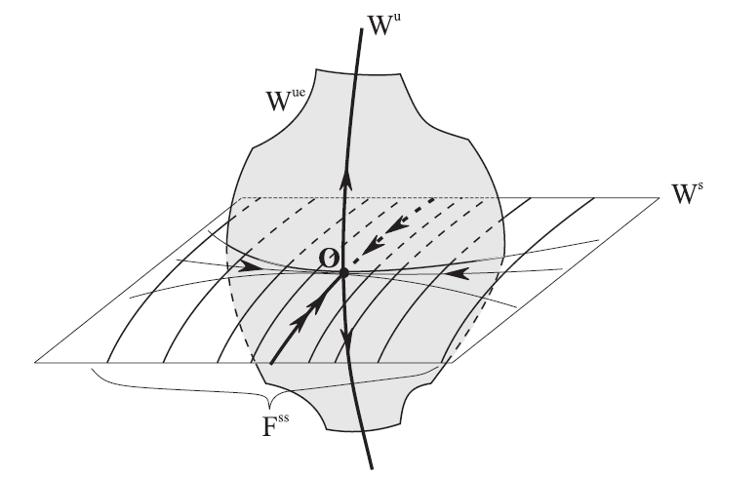,
height=6cm}} \caption{{\footnotesize A part of the strong stable
foliation $F^{ss}$ containing the strong stable manifold $W^{ss}$;
and a piece of one of the extended unstable manifolds containing
$W^u$ and being transversal to $W^{ss}$ at $O$.}} \label{fig04}
\end{figure}

According to \cite{book, HPS}, an important role in the global dynamics is played by local structures near the saddle fixed point. The first such object is an {\em extended unstable manifold }  $W^{ue}(O)$,
see figure~\ref{fig04}. By definition, it is a two-dimensional invariant manifold, tangent to the leading stable direction 
(those, corresponding to $\lambda_1$) at the saddle point and containing the unstable manifold $W^u(O)$. 
Unlike the previous ones, the extended unstable manifold is not uniquely defined and $C^{1 + \varepsilon}$--smooth. 
Locally, $W^{ue}_{loc}(O) = W^{ue}(O) \cap U_0$, and the
equation of $W^{ue}_{loc}(O)$ has the form $x_2 = \psi(x_1, y)$, where $\psi(0, y) \equiv 0$ and $\psi'_{x_1}(0, 0) = 0$.
Note that despite the fact that $W^{ue}(O)$ is non-unique, all of them have the same tangent plane at each point of $W^u(O)$.

Another essential structure is the {\em strong stable invariant foliation} on $W^s(O)$, see figure~\ref{fig04}.
In $W^s(O)$ there exists a one-dimensional strong stable invariant submanifold $W^{ss}(O)$, which is
$C^r$--smooth and touches at $O$ the strong stable direction. Moreover, manifold $W^s(O)$ is foliated near $O$ by
the leaves of an one-dimensional invariant foliation $F^{ss}$ which is $C^r$-smooth, unique and contains
$W^{ss}(O)$ as a leaf.

The existence of these structures implies that locally in the neighborhood $U_2$ the  dynamics can be locally decomposed into a dynamics on the two-dimensional invariant manifold $W^{ue}(O_2)$ and strong transversal contraction towards it along the leaves of the foliation $F^{ss}$. If such a decomposition spreads globally along the whole homoclinic (heteroclinic) cycle, it will imply the existence of a global two-dimensional invariant manifold. 
 In such a case, the existence of three-dimensional chaotic objects (discrete Lorenz attractors) is not possible. Thus, in order to have these attractors it is essential to prevent from existence of such global invariant manifolds. These conditions will be discussed in details in the next subsection.
 
If a fixed point experiences a transition from a saddle to a saddle-focus ($O$ or $O_1$), then for $\mu_2 > 0$, when it has real eigenvalues and becomes a saddle, the matrices  $A_1(\mu)$ and $A_2(\mu)$ are one-dimensional. But for $\mu \le 0$ their dimensions are different, namely $A_1(\mu)$ is two-dimensional and $A_2(\mu)$ zero-dimensional. It means that normal forms (\ref{t0norm}) for $\mu_2 < 0$ and $\mu_2 > 0$ can not be conjugated smoothly in a parametric family. So in order
to construct a smooth family of local maps, we need to use another form of the linear part. A normal form proposed in \cite{ArnGeom} for multiple eigenvalues plays well this role. Namely, the following formula explicitly defines $\mu_2$ and provides a smooth conjugation of families of matrices with real and complex-conjugated eigenvalues:
\begin{equation} \label{Alocal}
A_1(\mu) = \begin{pmatrix} \lambda & 1 \\ \mu_2 & \lambda
\end{pmatrix}.
\end{equation}
Indeed, when $\mu_2 > 0$ the eigenvalues of $A_1$ are real, $\lambda_{1,2} = \lambda \pm \sqrt{\mu_2}$ and for $\mu_2 < 0$ they are complex-conjugate
$\lambda_{1,2} = \lambda \pm i \sqrt{-\mu_2}$. To use this matrix in the first return map we need to find a suitable
form for its multiple iterations $A_1^k(\mu)$. Namely, if $\mu_2 > 0$:
\begin{equation} \label{Apowers+}
A_1^k(\mu) =  \begin{pmatrix} 
\displaystyle \frac{(\lambda + \sqrt{\mu_2})^k + (\lambda - \sqrt{\mu_2})^k}{2} & 
\displaystyle \frac{(\lambda + \sqrt{\mu_2})^k - (\lambda - \sqrt{\mu_2})^k}{2 \sqrt{\mu_2}} \\ 
\displaystyle  \sqrt{\mu_2}\frac{(\lambda + \sqrt{\mu_2})^k - (\lambda - \sqrt{\mu_2})^k}{2} & 
\displaystyle \frac{(\lambda + \sqrt{\mu_2})^k + (\lambda - \sqrt{\mu_2})^k}{2},
\end{pmatrix}
\end{equation}
if $\mu_2 < 0$
\begin{equation} \label{Apowers-}
A_1^k(\mu) =  \begin{pmatrix} 
\displaystyle \frac{(\lambda + i \sqrt{-\mu_2})^k + (\lambda - i \sqrt{-\mu_2})^k}{2} & 
\displaystyle \frac{(\lambda + i \sqrt{-\mu_2})^k - (\lambda - i \sqrt{-\mu_2})^k}{2i \sqrt{-\mu_2}} \\ 
\displaystyle  i\sqrt{-\mu_2}\frac{(\lambda + i\sqrt{-\mu_2})^k - (\lambda - i \sqrt{-\mu_2})^k}{2} & 
\displaystyle \frac{(\lambda + i \sqrt{-\mu_2})^k + (\lambda - i \sqrt{-\mu_2})^k}{2}
\end{pmatrix}
\end{equation}
 and if $\mu_2 = 0$:
\begin{equation} \label{Apowers0}
A_1^k(\mu) =  \begin{pmatrix} 
\lambda^k & k \lambda^{k - 1} \\ 
0 & \lambda^k
\end{pmatrix}.
\end{equation}
It is clear that $A_1^k(\mu)$ depends smoothly on $\mu_2$. Formulas (\ref{Apowers+})--(\ref{Apowers0}) can be
written in a more compact form. Define angle $\phi$ as follows:
\begin{equation} \label{phi}
\phi = \left\{
\begin{array}{l}
\displaystyle \arctanh \frac{\sqrt{\mu_2}}{\lambda}, \;\; \mbox{if} \;\; \mu_2 > 0 \\
0, \;\; \mbox{if} \;\; \mu_2 = 0 \\
\displaystyle -\arctan \frac{\sqrt{-\mu_2}}{\lambda}, \;\; \mbox{if} \;\; \mu_2 < 0. \\
\end{array}
\right.
\end{equation}
Then finally
\begin{equation} \label{ApowerC}
A^k(\mu) =  (\lambda^2 - \mu_2)^{k/2} \begin{pmatrix} 
C_{k}(\phi)   & S_{k}(\phi) \\ 
\mu_2 S_{k}(\phi) & C_{k}(\phi)
\end{pmatrix},
\end{equation}
where
\begin{equation} \label{Ckf}
C_k(\phi) = \left\{
\begin{array}{l}
\cosh k \phi, \;\; \mbox{if} \;\; \mu_2 \ge 0 \\
\cos k \phi, \;\; \mbox{if} \;\; \mu_2 < 0. \\
\end{array}
\right.
\end{equation}
and
\begin{equation} \label{Skf}
S_k(\phi) = \left\{
\begin{array}{l}
\displaystyle \frac{\sinh k \phi}{\sqrt{\mu_2}}, \;\; \mbox{if} \;\; \mu_2 > 0 \\
\displaystyle \frac{k}{\lambda}, \;\; \mbox{if} \;\; \mu_2 = 0 \\
\displaystyle -\frac{\sin k \phi}{\sqrt{-\mu_2}}, \;\; \mbox{if} \;\; \mu_2 < 0. \\
\end{array}
\right.
\end{equation}
are smooth functions of $\mu_2$.
Now, multiple iterations near $O$, 
$T_{0}^k : U_{0} \to U_{0}$ (or near $O_1$), for any $k$ can be also calculated in a form close to linear. 
For all small $\mu$ iterations $T_{0}^k: (x_0, y_0) \to (x_k, y_k)$ can be represented as:
\begin{equation}\label{eq:T0kk1}
\begin{array}{l}
  (x_{1k}, x_{2k})^\top = A_1^k(\mu) (x_{10}, x_{20})^\top +
            \hat\lambda^{k} \xi_{1k}(x_0, y_k, \mu)\\
  y_0  = \gamma^{-k} y_k + \hat \gamma^{-k} \xi_{2k}(x_0, y_k, \mu).
\end{array}
\end{equation}
Here also, $0 < \hat \lambda < \lambda$, $\hat\gamma > \gamma$, 
 functions $\xi_{k}$ are uniformly bounded together with their derivatives up to order $(r - 2)$
 and their higher-order derivatives tend to zero.


\subsection{Properties of the global maps}

Consider the homoclinic Case I.
Assume that for
$\mu = 0$, the coordinates of the homoclinic points are $M^-(0, 0, y^-) \in W^u_{loc}(O)$, $M^+(x_1^+, x_2^+, 0)\in W^s_{loc}(O)$.
Then the global map for all small $\mu$ can be written as a Taylor expansion near point $M^-$:
\begin{equation}\label{eq:T1}
  \begin{array}{l}
    \bar x_{10} - x_1^+ = a_{11} x_{1k} +a_{12}x_{2k} + b_1 (y_k - y^-) + O(\|x_k\|^2 + \|x_k\| \cdot |y_k - y^-| + (y_k - y^-)^2), \\
    \bar x_{20} - x_2^+ = a_{21} x_{1k} +a_{22}x_{2k} + b_2 (y_k - y^-) + O(\|x_k\|^2 + \|x_k\| \cdot |y_k - y^-| + (y_k - y^-)^2), \\

    \bar y_0 = y^+ + c_1 x_{1k} + c_2 x_{2k} + d (y_k - y^-)^2 +  O(\|x_k\|^2 + \|x_k\| \cdot |y_k - y^-| + (y_k - y^-)^3),
  \end{array}
\end{equation}
where all the coefficients smoothly depend on parameters. We note also that $d \neq 0$ since $W^u(O)$ and $W^s(O)$ 
have a quadratic tangency for $\mu = 0$, i.e. $y^+(0) = 0$. With that we define the first control parameter to be
$\mu_1 \equiv y^+(\mu)$. Map $T_{1}$ is a diffeomorphism, so that
\begin{equation}\label{eq:J1}
J_{T1} = \mbox{det}\; \begin{pmatrix}
a_{11} & a_{12} & b_1 \\
a_{21} & a_{22} & b_2 \\
c_1 & c_2 & 0
\end{pmatrix}
\; \neq 0
\end{equation}
In particular, this means that $b_1^2 + b_2^2 \neq 0$ and $c_1^2 + c_2^2 \neq 0$ for $\mu = 0$.

For small $\mu_2 > 0$ the fixed point $O$ becomes a saddle point. Therefore, 
additional degeneracies may occur with the homoclinic orbit. Namely, the homoclinic tangency can become non-simple. Next we give the definition of it and formulate the non-degeneracy condition in coordinates.

Denote as $P^{ue}(M^-)$ a tangent plane to to $W^{ue}(O)$ at the point $M^-$. 
The homoclinic tangency of $W^u(O)$ and $W^s(O)$ is called {\em simple} if image $T_{1}(P^{ue}(M^-))$
intersects transversely (i.e. generically) the  leaf $F^{ss}(M^+)$ of invariant foliation $F^{ss}$, 
coming through point $M^+$. If this condition is not fulfilled, such a quadratic tangency is called {\em non-simple}. There are two possibilities of codimension-one non-simple homoclinic orbits: 1) inclination flip, when the surface $T_{1}(P^{ue}(M^-))$ is tangent to surface $W^s(O)$, but $W^u(O)$ is not tangent to the leaf $F^{ss}(M^+)$, and 2) orbit flip, when $T_{1}(P^{ue}(M^-))$ is transverse to $W^s(O)$, but $W^u(O)$ is tangent to the leaf $F^{ss}(M^+)$. It is clear that in both cases further forward iterations of $W^{ue}(O)$ will not tend to be tangent to $W^{ue}_{loc}(O)$, hence a global two-dimensional invariant manifold will not exist.

The notion of a simple quadratic homoclinic tangency was introduced in \cite{GST93c}, and it is a variant of the so-called
quasi-transversal intersection \cite{NPT}. A simple homoclinic tangency of arbitrary order
was defined in \cite{GLi12}. Bifurcations of non-simple homoclinic tangencies were
studied first in \cite{Tat01}, where it was called ``generalized homoclinic tangency'', and also \cite{GGT07}. Birth of discrete Lorenz attractors in bifurcations of non-simple homoclinic and heteroclinic orbits were studied in \cite{GOT14, O22}.

Now consider small $\mu_2 > 0$ and a homoclinic tangency (i.e. $\mu_1 = 0$), for which $O$ is a saddle with stable multipliers $\lambda_{1,2} = \lambda \pm \sqrt{\mu_2}$ and
unstable one $\gamma$. In the notation of section~\ref{sec:def}, this means that 
$\dim W^{ss}(O) = \dim (W^{s}(O) \backslash W^{ss}(O)) = \dim W^u(O) = 1$, $\dim W^{ue}(O) = 2$ and 
invariant foliation $F^{ss}$ consists of one-dimensional leaves.
The stable eigenvectors at $O$ are: leading $l^s = (1, \sqrt{\mu_2}, 0)^\top$ 
and non-leading $l^{ss} = (-1, \sqrt{\mu_2}, 0)^\top$. The leaves of strong stable invariant foliation $F^{ss}$
are straight lines parallel to $l^{ss}$ and tangent plane $P^{ue}$ to extended unstable manifold $W^{ue}(O)$ has equation
$x_2 = x_1\sqrt{\mu_2} $. The homoclinic tangency will be simple if the following two conditions hold: tangent vector $l^u$ to
$W^u(O)$ is not parallel to $l^{ss}$ and surface $T_1(P^{ue})$ intersect $W^s(O)$ transversely. The latter will be fulfilled if
there is a tangent vector to $P^{ue}$, whose image under $T_1$ has a non-zero $y$--component.

Tangent to $W^u(O)$ vector $l^u$ has coordinates $(0, 0, 1)^\top$. Its image under $T_1$ is obtained by means of the differential of (\ref{eq:T1})
evaluated at $M^-$. It gives vector $DT_1(M^-)l^u = (b_1, b_2, 0)^\top$. This vector will be parallel to $l^{ss}$ (orbit flip) and the tangency will be non-simple 
if $b_1 \sqrt{\mu_2} - b_2 = 0$. Thus, condition $b_2 \neq 0$ guarantees that orbit flip does not happen for sufficiently small $\mu_2$.
Now take any tangent to $P^{ue}$ vector $(x_1, x_1 \sqrt{\mu_2}, y)$ with non-zero $x_1$ (recall that the image of vector $(0, 0, y)$ lies in $W^s$ at $\mu = 0$) and apply $DT_1(M^-)$ to it.
The $y$--component of the resulting vector will be $(c_1 + c_2 \sqrt{\mu_2}) x_1$. If this expression is zero for some parameter value, then surfaces $T_1 P^{ue}$ is tangent to $W^s(O)$, and the tangency is also non-simple (inclination flip). This additional degeneracy will not happen  
for sufficiently small $\mu_2$ if $c_1 \neq 0$. Thus we have the following non-degeneracy conditions:
\begin{equation}\label{Simple_hom}
b_2(0) \neq 0, \;\; c_1(0) \neq 0.
\end{equation}

Thus, the global map $T_1$ in Case I is given by equation (\ref{eq:T1}) with non-degeneracy conditions (\ref{eq:J1}) and (\ref{Simple_hom}).

Now consider the heteroclinic Case II. 
Let for $\mu = 0$ the previously selected heteroclinic points have the following coordinates:
$M_1^-(0, 0, y^-)\in W^u_{loc}(O_1)$, 
$M_1^+(x_1^+, x_2^+, 0)\in W^s_{loc}(O_1)$,
$M_2^-(0, 0, v^-)\in W^u_{loc}(O_2)$,
$M_2^+(u_1^+, u_2^+, 0)\in W^s_{loc}(O_2)$.
The heteroclinic cycle is composed of a fixed point $O_1$ that changes from saddle to saddle-focus when $\mu_2$ changes sign, and a saddle point $O_2$. 

In the Case IIa heteroclinic orbit $\Gamma_{12}$ is transversal and heteroclinic orbit $\Gamma_{21}$ consists of quadratic tangencies. Then, the global map $T_{12}$ can be written for all small $\mu$ as a Taylor expansion near point $M_1^-$:
\begin{equation}\label{eq:T12a}
  \begin{array}{l}
    u_{10} - u_1^+ = a_{11} x_{1k} +a_{12}x_{2k} + b_1 (y_k - y^-) + O(\|x_k\|^2 + \|x_k\| \cdot |y_k - y^-| + (y_k - y^-)^2), \\
    u_{20} - u_2^+ = a_{21} x_{1k} +a_{22}x_{2k} + b_2 (y_k - y^-) + O(\|x_k\|^2 + \|x_k\| \cdot |y_k - y^-| + (y_k - y^-)^2), \\

    v_0 = c_1 x_{1k} + c_2 x_{2k} + d_1 (y_k - y^-) +  O(\|x_k\|^2 + \|x_k\| \cdot |y_k - y^-| + (y_k - y^-)^2).
  \end{array}
\end{equation}
The global map $T_{21}$ is written as a Taylor expansion near point $M_2^-$:
\begin{equation}\label{eq:T21a}
  \begin{array}{l}
    \bar x_{10} - x_1^+ = a_{31} u_{1j} +a_{32}u_{2j} + b_3 (v_j - v^-) + O(\|u_j\|^2 + \|u_j\| \cdot |v_j - v^-| + (v_j - v^-)^2), \\
    \bar x_{20} - x_2^+ = a_{41} u_{1j} +a_{42}u_{2j} + b_4 (v_j - v^-) + O(\|u_j\|^2 + \|u_j\| \cdot |v_j - v^-| + (v_j - v^-)^2), \\

    \bar y_0 = \mu_1 + c_3 u_{1j} + c_4 u_{2j} + d_2 (v_j - v^-)^2 +  O(\|u_j\|^2 + \|u_j\| \cdot |v_j - v^-| + (v_j - v^-)^3).
  \end{array}
\end{equation}
The non-degeneracy conditions on Jacobians have the form:
\begin{equation}\label{eq:J12a}
J_{12} = \mbox{det}\; \begin{pmatrix}
a_{11} & a_{12} & b_1 \\
a_{21} & a_{22} & b_2 \\
c_1 & c_2 & d_1
\end{pmatrix}
\; \neq 0,
\end{equation}
\begin{equation}\label{eq:J21a}
J_{21} = \mbox{det}\; \begin{pmatrix}
a_{31} & a_{32} & b_3 \\
a_{41} & a_{42} & b_4 \\
c_3 & c_4 & 0
\end{pmatrix}
\; \neq 0
\end{equation}
and also
\begin{equation}\label{d1d2a}
    d_1 \neq 0, \; d_2 \neq 0.
\end{equation} 

Other possible degeneracies that should be prevented, are related to the existence of non-simple heteroclinic orbits. For quadratic heteroclinic tangencies the condition to be non-simple is defined in a similar way as for homoclinic ones -- we consider the extended unstable manifold $W^{ue}(O_2)$ and its image under $T_{21}$. If this image is transverse to the leaf $F^{ss}(M_1^+)$, the tangency is called simple. Similarly, we obtain two possibilities of non-simple tangencies -- the orbit flip and the inclination flip. For transversal heteroclinic orbits, it is also possible to have an orbit flip, when image $T_{12}W^{ue}(O_1)$ is tangent to the leaf $F^{ss}(M_2^+)$. Such bifurcation was  studied first in \cite{O22}.

For the orbit $\Gamma_{12}$, the extended unstable manifold $W^{ue}(O_1)$ is tangent to $x_1$ direction as $\mu_2 \to 0$, and the leaf of invariant foliation $F^{ss}$ at point $M^+_2$ is $\{v = 0, \; u_1 = u_1^+ \}$, the tangent vector to this leaf is $l_{ss} = (0, 1, 0)^\top$. The image $T_{12} P^{ue}(M_1^-)$ has two linearly independent tangent vectors $l_1 = (a_{11}, a_{21}, c_1)^\top$ and $l_2 = (b_1, b_2, d_1)^\top$. The heteroclinic intersection will be simple if vectors $l_1$, $l_2$ and $l_{ss}$ do not lie in the same plane, that is:
\begin{equation}\label{noOrbitFlip12a}
   A_{11}( 0) = \left. \left(a_{11} - \frac{b_1 c_1}{d_1}\right) \right|_{\mu = 0} \neq 0.
\end{equation}
%
%

For the orbit $\Gamma_{21}$  the extended unstable manifold $W^{ue}(O_2)$ is tangent to $u_1$ direction, and the leaves of invariant foliation $F^{ss}$ in $U_1$ tend to be parallel to $x_1$ as $\mu_2 \to 0$, the tangent vector to this leaf is $l_{ss} = (1, 0, 0)^\top$. The image $T_{21} P^{ue}(M_2^-)$ has two linearly independent tangent vectors $l_1 = (a_{31}, a_{41}, c_3)^\top$ and $l_2 = (b_3, b_4, 0)^\top$. The heteroclinic intersection will be simple if vectors $l_1$, $l_2$ and $l_{ss}$ do not lie in the same plane,
that is:
\begin{equation}\label{Simple21a}
   b_4(0) c_3(0) \neq 0.
\end{equation}

To summarize, in the Case IIa the global maps are given by formulas (\ref{eq:T12a}) and (\ref{eq:T21a}) with non-degeneracy conditions (\ref{eq:J12a})--(\ref{Simple21a}).

Analogously, we construct the global maps in the Case IIb, when orbit $\Gamma_{12}$ consists of quadratic tangencies and $\Gamma_{21}$ is transversal. The global maps $T_{12}$ and $T_{21}$ have the following form:
\begin{equation}\label{eq:T12b}
  \begin{array}{l}
    u_{10} - u_1^+ = a_{11} x_{1k} +a_{12}x_{2k} + b_1 (y_k - y^-) + O(\|x_k\|^2 + \|x_k\| \cdot |y_k - y^-| + (y_k - y^-)^2), \\
    u_{20} - u_2^+ = a_{21} x_{1k} +a_{22}x_{2k} + b_2 (y_k - y^-) + O(\|x_k\|^2 + \|x_k\| \cdot |y_k - y^-| + (y_k - y^-)^2), \\

    v_0 = \mu_1 + c_1 x_{1k} + c_2 x_{2k} + d_1 (y_k - y^-)^2 +  O(\|x_k\|^2 + \|x_k\| \cdot |y_k - y^-| + (y_k - y^-)^3).
  \end{array}
\end{equation}
and
\begin{equation}\label{eq:T21b}
  \begin{array}{l}
    \bar x_{10} - x_1^+ = a_{31} u_{1j} +a_{32}u_{2j} + b_3 (v_j - v^-) + O(\|u_j\|^2 + \|u_j\| \cdot |v_j - v^-| + (v_j - v^-)^2), \\
    \bar x_{20} - x_2^+ = a_{41} u_{1j} +a_{42}u_{2j} + b_4 (v_j - v^-) + O(\|u_j\|^2 + \|u_j\| \cdot |v_j - v^-| + (v_j - v^-)^2), \\

    \bar y_0 =  c_3 u_{1j} + c_4 u_{2j} + d_2 (v_j - v^-) +  O(\|u_j\|^2 + \|u_j\| \cdot |v_j - v^-| + (v_j - v^-)^2),
  \end{array}
\end{equation}
and the non-degeneracy conditions are the non-zero Jacobians:
\begin{equation}\label{eq:J12b}
J_{12} = \mbox{det}\; \begin{pmatrix}
a_{11} & a_{12} & b_1 \\
a_{21} & a_{22} & b_2 \\
c_1 & c_2 & 0
\end{pmatrix}
\; \neq 0,
\end{equation}
and
\begin{equation}\label{eq:J21b}
J_{21} = \mbox{det}\; \begin{pmatrix}
a_{31} & a_{32} & b_3 \\
a_{41} & a_{42} & b_4 \\
c_3 & c_4 & d_2
\end{pmatrix}
\; \neq 0
\end{equation}
the conditions of quadratic tangency and transverse orbits:
\begin{equation}\label{d1d2b}
    d_1 \neq 0, \; d_2 \neq 0,
\end{equation}
and absence of non-simple orbits:
\begin{equation}\label{Simple12b}
   b_1(0) c_1(0) \neq 0.
\end{equation}
and
\begin{equation}\label{noOrbitFlip12b}
   A_{41}( 0) = \left. \left(a_{41} - \frac{b_4 c_3}{d_2}\right) \right|_{\mu = 0} \neq 0.
\end{equation}

\begin{lm}[Rescaling lemma]\label{th4}
Let $f_{\mu_1, \mu_2, \mu_3}$ be the family under consideration. Then, in the parameter
space $(\mu_1, \mu_2, \mu_3)$ in Case I there exist infinitely many regions $\Delta_k$
ac\-cu\-mu\-la\-ting to the origin as $k\to \infty$, and in Case II there exist infinitely many regions $\Delta_{kj}$
ac\-cu\-mu\-la\-ting to the origin as $k, j\to \infty$, such that the first return map  in appropriate rescaled
coordinates and parameters is asymptotically $C^{r - 1}$-close to the  map (\ref{H3D}).
%
%
with parameters

{\rm 1)} In Case I:
%
\begin{equation}
\label{M123I}
\begin{array}{l}
M_1 = -d^2 \gamma^{2k}
     (\mu_1 + \ldots)\\

M_2 = L(\mu)b_2  \lambda^k \gamma^k \sin(k \phi + \psi) + \ldots  \\
B = J_{T1}(\lambda^2\gamma)^k + \ldots
\end{array}
\end{equation}

{\rm 2)}  In Case IIa:
%
\begin{equation}
\label{M123IIa}
\begin{array}{l}
M_1 = -d_{1}^2 d_{2} \gamma_1^{2k} \gamma_2^{2j}
     (\mu_1 + \ldots)\\
M_2 = L(\mu) b_4 d_1  (\lambda\gamma_1)^k(\nu_1\gamma_2)^j \sin(k \phi + \psi) + \ldots \\
B = J_{12}J_{21}(\lambda^2\gamma_1)^k(\nu_1\nu_2\gamma_2)^j + \ldots
\end{array}
\end{equation}
%

{\rm 3)}  In Case IIb:
%
\begin{equation}
\label{M123IIb}
\begin{array}{l}
M_1 = -d_{1}^2 d_{2} \gamma_1^{2k} \gamma_2^{2j}
     (\mu_1 + \ldots)\\

M_2 = L(\mu) b_1 d_2 A_{41}  (\lambda\gamma_1)^k(\nu_1\gamma_2)^j \sin(k \phi + \psi) + \ldots \\
B = J_{12}J_{21}(\lambda^2\gamma_1)^k(\nu_1\nu_2\gamma_2)^j + \ldots
\end{array}
\end{equation}
where $L(\mu)$ is some bounded from zero function, the parameters are varied such that $\sin(k \phi + \psi)$ is kept small, and dots stand for terms that tend to zero as  $k,j \to \infty$.

\end{lm}

The rescaling normal form for the first 
return map $T_k$ or $T_{kj}$ is the 3D H\'enon map (\ref{H3D}),
in  both Cases~$\rm I$ and~$\rm II$. The statement of the main Theorem~\ref{thm:main} immediately follows from this lemma: indeed, in the H\'enon map (\ref{H3D}) there exist  parameter domains $\Delta^+$ and $\Delta^-$ (for $B > 0$ and $B < 0$ respectively), at which the map possesses a discrete Lorenz attractor. Then, in the Case I for every sufficiently large $k$ they correspond to subsets $\delta_k$ that converge to the origin in the space of original parameters $(\mu_1, \mu_2, \mu_3)$ as $k \to \infty$. In the heteroclinic Case II, for every sufficiently large $k, j$ such that the product of Jacobians of local maps $(\lambda^2\gamma_1)^k(\nu_1\nu_2\gamma_2)^j$ is kept finite, domains $\Delta^\pm$ correspond to subsets $\delta_{kj}$ that also converge to the origin as $k, j \to \infty$.


\section{Proof of the rescaling  Lemma.}\label{sec:th3proof}

\subsection{Homoclinic Case I}

Using formulas (\ref{eq:T0kk1}), (\ref{eq:T1}), we obtain the following expression for the first
return map $ T_{k} \equiv T_{1} T_{0}^k: \sigma_0^k \subset \Pi^+ \to\Pi^+$
\begin{equation}\label{eq:Tk}
  \begin{array}{l}
    \bar x_1 - x_1^+ = a_{11} \lambda^k (x_1 C_k(\phi) + x_2 S_k(\phi)) + a_{11}\hat\lambda^{k} \xi_{1k}^1(x, y, \mu) + \\ 
    \qquad \qquad + a_{12} \lambda^k  (x_1 \mu_2 S_k(\phi) + x_2 C_k(\phi)) + a_{12}\hat\lambda^{k} \xi_{1k}^2(x, y, \mu) + \\
    \qquad \qquad + b_1 (y - y^-) + O(\lambda^{2k} \|x\|^2 + \lambda^k \|x\| \cdot |y - y^-| + (y - y^-)^2), \\
    
    \bar x_2 - x_2^+ = a_{21} \lambda^k  (x_1 C_k(\phi) + x_2 S_k(\phi)) + a_{21}\hat\lambda^{k} \xi_{1k}^1(x, y, \mu) + \\   
    \qquad \qquad + a_{22} \lambda^k  (x_1 \mu_2 S_k(\phi) + x_2 C_k(\phi)) + a_{22}\hat\lambda^{k} \xi_{1k}^2(x, y, \mu) + \\   
    \qquad \qquad + b_2 (y - y^-) + O(\lambda^{2k}\|x\|^2 + \lambda^k \|x\| \cdot |y - y^-| + (y - y^-)^2), \\

   \gamma^{-k} (\bar y + \left( \hat\gamma / \gamma \right)^{-k} \xi_{2k}(\bar x, \bar y, \mu)) = \mu_1 + 
    c_1 \lambda^k  (x_1 C_k(\phi) + x_2 S_k(\phi)) + \\ 
    \qquad \qquad + c_1 \hat\lambda^{k} \xi_{1k}^1(x, y, \mu)) + c_2 \lambda^k  (x_1 \mu_2 S_k(\phi) + x_2 C_k(\phi)) + \\
    \qquad \qquad + c_{2}\hat\lambda^{k} \xi_{1k}^2(x, y, \mu) + d(y - y^-)^2 + \\
    \qquad \qquad + O(\lambda^{2k}\|x\|^2 + \lambda^k \|x\| \cdot |y - y^-| + (y - y^-)^3), \\
  \end{array}
\end{equation}

Make a coordinate shift
$x_{1new} = x_1 - x_1^+ + \psi_{k}^1$, $x_{2new} = x_2 - x_2^+ + \psi_{k}^2$, $y_{new} = y - y^- + \psi_{k}^3$,
where $\psi_{k}^{1,2,3} = O(\lambda^k)$. With that, the nonlinearity function in the left part of the third equation
of (\ref{eq:Tk}) can be expressed as a Taylor expansion
$\xi_{2k}(\bar x_1 + x_1^+ - \psi_{k}^1, \bar x_2 + x_2^+ - \psi_{k}^2, \bar y + y^- - \psi_{k}^3, \mu)) = 
\xi_{k}^0 + \xi_{k}^1 \bar y + \xi_{k}^2(\bar x, \bar y) + \xi_{k}^3(\bar y)$, where coefficients  
$\xi_{k}^0$, $\xi_{k}^1$ are uniformly bounded in $k$ for all small $\mu$ and
$\xi_{k}^2 (\bar x, \bar y) = O(\bar x)$, $\xi_{k}^3(\bar y) = O(\bar y^2)$. 
We select constants $\psi_{k}^1$, $\psi_{k}^2$, $\psi_{k}^3$ in such a way
that all constant terms in the first two equations and the linear in $y$ term
in the last equation of (\ref{eq:Tk}) vanish. Next, we substitute the right hand sides of first two equations to term
$O(\bar x)$ in the last one. All coefficients then get additions of order $\hat\gamma^{-k}$ including the appearance
of one linear in $y$ in the third equation. The
latter can be again made zero by an additional shift of $(x_1, x_2, y)$ coordinates by the value of order $\hat\gamma^{-k}$.
The system is rewritten as:
\begin{equation}\label{eq:3}
    \begin{array}{l}
    \bar x_1 = \lambda^k  (a_{11} C_k(\phi) + a_{12} \mu_2 S_k(\phi)) x_1 + 
    \lambda^k  (a_{11} S_k(\phi) + a_{12} C_k(\phi)) x_2 + \\
    \qquad \qquad + b_1 y + \hat\lambda^k O(\|x\|^2) + \lambda^k O(\|x\| \cdot |y|) + O(y^2), \\
    
    \bar x_2 = \lambda^k  (a_{21} C_k(\phi) + a_{22} \mu_2 S_k(\phi)) x_1 + 
    \lambda^k (a_{21} S_k(\phi) + a_{22}C_k(\phi)) x_2 + \\
    \qquad \qquad + b_2 y + \hat\lambda^k O(\|x\|^2) + \lambda^k O(\|x\| \cdot |y|) + O(y^2), \\

    \gamma^{-k} (1 + q_{k})\bar y + \hat\gamma^{-k} O(\bar y^2) = M^1 +
    \lambda^k  (c_{1} C_k(\phi) + c_{2} \mu_2 S_k(\phi)) x_1 + \\    
    \qquad \qquad +\lambda^k  (c_{1} S_k(\phi) + c_{2} C_k(\phi)) x_2 + d y^2 + \\
    \qquad \qquad + \hat\lambda^k O(\|x\|^2) + \lambda^k O(\|x\| \cdot |y|) + O(y^3), \\
  \end{array}
\end{equation}
where $q_{k} = O\left(\left( \hat\gamma / \gamma \right)^{-k} \right)$ and the following expression is valid for $M^1$:

\begin{equation} \label{MfFeps}
M^1 = \mu_1 + \lambda^k c_1 S_k(\phi) x_2^+ - \gamma^{-k}y^- + 
O(\hat\gamma^{-k} + \lambda^k). 
\end{equation}

Next, as $b_2 \neq 0$, we make the following change of coordinates in order to eliminate terms that depend on $y$ only in the first equation:
$\displaystyle x_{1new} = x_1 - \frac{b_1}{b_2} x_2 + O(\|x\|^2)$.
System (\ref{eq:3}) takes the following form:
\begin{equation} \label{eq:4}
    \begin{array}{l}
    \bar x_1 = \lambda^k  A_{11k}(\phi) x_1 + \lambda^k  A_{12k}(\phi) x_2 
 + \hat\lambda^k O(\|x\|^2) + \lambda^k O(\|x\| \cdot |y|) , \\
    
    \bar x_2 = \lambda^k  A_{21k}(\phi) x_1 + \lambda^k  A_{22k}(\phi) x_2 +
    b_2 y + \hat\lambda^k O(\|x\|^2) + \\ 
    \qquad \qquad   + \lambda^k O(\|x\| \cdot |y|) + O(y^2), \\

    \gamma^{-k} (1 + q_{k})\bar y + \hat\gamma^{-k} O(\bar y^2) = M^1 +
    \lambda^k  C_{1k}(\phi) x_1 + \lambda^k  C_{2k}(\phi) x_2 + \\ 
    \qquad \qquad + d y^2 + \hat\lambda^k O(\|x\|^2) + \lambda^k O(\|x\| \cdot |y|) + O(y^3), \\
  \end{array}
\end{equation}
where
\begin{equation}
\begin{array}{l}
\displaystyle A_{11k}(\phi) = a_{11} C_k(\phi) + a_{12} \mu_2 S_k(\phi) - \frac{b_1}{b_2} (a_{21} C_k(\phi) + a_{22}\mu_2 S_k(\phi) )\\
\displaystyle A_{12k}(\phi) = a_{11} \left( \frac{b_1}{b_2} C_k(\phi) + S_k(\phi) \right) + a_{12} \left( \frac{b_1}{b_2}\mu_2 S_k(\phi)
+ C_k(\phi) \right) - \\
\displaystyle \qquad - a_{21} \frac{b_1}{b_2} \left( \frac{b_1}{b_2} C_k(\phi) + S_k(\phi) \right) 
       - a_{22} \frac{b_1}{b_2} \left( \frac{b_1}{b_2} \mu_2 S_k(\phi) + C_k(\phi) \right)\\
       
\displaystyle A_{21k}(\phi) = a_{21} C_k(\phi) + a_{22} \mu_2 S_k(\phi)\\
\displaystyle A_{22k}(\phi) = a_{21} \left( \frac{b_1}{b_2} C_k(\phi) + S_k(\phi) \right) +
         a_{22} \left( \frac{b_1}{b_2} \mu_2 S_k(\phi) + C_k(\phi) \right) \\
\displaystyle C_{1k}(\phi) = c_1 C_k(\phi) + c_2 \mu_2 S_k(\phi) \\
\displaystyle C_{2k}(\phi) = c_1 \left( \frac{b_1}{b_2} C_k(\phi) + S_k(\phi) \right) + 
          c_{2} \left( \frac{b_1}{b_2} \mu_2 S_k(\phi) + C_k(\phi) \right)
\end{array}
\end{equation}

Now consider $C_{2k}(\phi)$ for $\mu_2 < 0$. It can be represented in the form $C_{2k}(\phi) = L_1(\mu) \cos k \phi + L_2(\mu) \sin k \phi = L(\mu) \sin(k \phi + \psi)$. We will vary $\mu_2$ near those values, for which $C_{2k}(\phi)$ vanishes (such values are dense on any interval from $(0, \pi)$), we can make this coefficient asymptotically small for $k \to \infty$. Note that for such values of $\mu_2$ expressions $\displaystyle \frac{b_1}{b_2} C_k(\phi) + S_k(\phi)$ and $\displaystyle \frac{b_1}{b_2} \mu_2 S_k(\phi) + C_k(\phi)$ are non-zeros, because otherwise coefficients $A_{12k}(\phi)$, $A_{21k}(\phi)$ and $C_{2k}(\phi)$ will vanish simultaneously, and the Jacobian of the first return map $T_k$ will become zero, which is impossible because it is a diffeomorphism. Then, it is easy to check that 
\begin{equation}\label{A12nonzeroI}
 \displaystyle   A_{12k}(\phi) = \frac{J_{T1}}{b_2 c_1}\left( \frac{b_1}{b_2} \mu_2 S_k(\phi) + C_k(\phi) \right) + O(C_{2k}(\phi)) \neq 0.
\end{equation}
Now we make the following scaling of coordinates: $x_1 = \alpha_1 X_1, \; x_2 = \alpha_2 X_2, \; y = \beta Y$, with
\begin{equation}\label{scale1}
\displaystyle \alpha_1 = - \frac{A_{12k} b_2 \gamma^{-k} \lambda^k}{d }(1 + q_k), \;\; 
\alpha_2 = -\frac{b_2 \gamma^{-k}}{d}(1 + q_k), \;\; \beta = -\frac{\gamma^{-k}}{d } (1 + q_k).
\end{equation}
%
The system becomes the following:
\begin{equation} \label{3DHenon1}
\begin{array}{l}
\bar X_1 = X_2 + O(\lambda^k) \\
\bar X_2 = Y + O(k \lambda^k) \\
\bar Y = M_1 + J_k X_1 + M_2 X_2 - Y^2 + O(\lambda^k),
\end{array}
\end{equation}
where $\displaystyle M_1 = -d \gamma^{2k} M^1$, $M_2 = \lambda^k  \gamma^k b_2 C_{2k}(\phi)$ and $J_k = (\lambda^2 \gamma)^k J_1  + O(k \lambda^k)$ coincides with the 
Jacobian of the first return map $T_k$ up to asymptotically small terms. For arbitrary large $k$ multiplier $\gamma^{2k}$ unboundedly
grows, so small variations of $\mu_1$ will provide arbitrary finite values for $M_1$. Then, as $|\lambda^2 \gamma| = 1 + \mu_3$ at the
bifurcation moment, the product $\lambda^k \gamma^k \sim \lambda^{-k}$ grows unboundedly as $k \to +\infty$. 
Varying parameter $\mu_2$ such that 
coefficient $C_{2k}(\phi)$ is kept small, we can make parameter $M_2$ to take arbitrary finite values. The third parameter 
$J_k$ can take arbitrary finite positive values when $\mu_3$ gets small perturbations.
System (\ref{3DHenon1}) is asymptotically close to (\ref{H3D}) as $k,j \to \infty$.

\begin{rem}
When $\mu_2 > 0$, scaling (\ref{scale1}) does not bring map $T_k$ to the 3D Henon form 
because in that case coefficient $C_{2k}(\phi)$ becomes growing in $k$ for all small $\mu_2$, 
so that parameter $M_2$ will be asymptotically large. In order to understand the dynamics in this case, 
in map (\ref{eq:4}) we perform an additional change of coordinates, which is valid because $C_{2k}(\phi) \neq 0$: 
\begin{equation}
\displaystyle x_{2new} = x_2 + \frac{C_{1k}(\phi)}{C_{2k}(\phi)} x_1 + 
O\left((\hat\lambda / \lambda)^{k} \right) O(\|x\|^2)
\end{equation}
in order to eliminate in the last equation the linear in $x_1$ term together with higher-order terms depending on $x$ only.
The system rewrites as:
\begin{equation} \label{eq:5}
    \begin{array}{l}
    \displaystyle \bar x_1 = \lambda^k  \frac{J_1}{b_2 C_{2k}(\phi)} x_1 + \lambda^k  A_{12k}(\phi) x_2 + 
 \hat\lambda^k O(\|x\|^2) + \lambda^k O(\|x\| \cdot |y|) , \\
    
    \bar x_2 = O(k \lambda^k) O(\|x\|) +  b_2 y + \hat\lambda^k O(\|x\|^2) 
     + \lambda^k O(\|x\| \cdot |y|) + O(y^2), \\

     (1 + q_{k})\bar y + \hat\gamma^{-k} O(\bar y^2) = \gamma^{k} M^1 +
     \lambda^{k}\gamma^{k} C_{2k}(\phi) x_2 + d\gamma^{k} y^2  + \\
     \qquad \qquad + \lambda^k\gamma^{k} O(\|x\| \cdot |y|) + \gamma^{k}O(y^3). \\
  \end{array}
\end{equation}

Finally, we scale the variables as follows: $x_1 = \alpha_1 X_1, \; x_2 = \alpha_2 X_2, \; y = \beta Y$, where
\begin{equation}\label{scale2}
\displaystyle \alpha_1 = \alpha_2 = -\frac{b_2 \gamma^{-k}}{d }(1 + q_k), \;\; \beta = -\frac{\gamma^{-k}}{d }(1 + q_k).
\end{equation}
Map (\ref{eq:5}) rewrites as follows:
\begin{equation} \label{2DHenon1}
\begin{array}{l}
\bar X_1 = O(k \lambda^k) \\
\bar X_2 = Y + O(k \lambda^k) \\
\bar Y = M_1 + M_2 X_2 - Y^2 + O(\lambda^k),
\end{array}
\end{equation}
where $\displaystyle M_1 = -d \gamma^{2k} M^1$, $M_2 = \lambda^k  \gamma^k b_2 C_{2k}(\phi)$. 
From formula (\ref{2DHenon1}) it is clear that the first return map has an invariant manifold of the form 
$X_1 = k \lambda^k \psi(X_2, Y)$ on each $\sigma$ and the restriction on it of system (\ref{2DHenon1}) is the well-known two-dimensional Henon map
with an arbitrary large Jacobian $M_2$ for large $k$, so it expands volumes and does not have attractors. 
\end{rem}


\subsection{Heteroclinic Case IIa.}\label{sect:caseIIa}

The first return map is a composition of two local and two global maps, defined above. We substitute  local maps  (\ref{eq:T0kk1}) and (\ref{eq:T20j}) to  global maps 
(\ref{eq:T12a}) and 
(\ref{eq:T21a}). Note that for simplicity we omit indices $0$, $k$ and $j$ at the variables:
\begin{equation}\label{eq:Tkj12a}
    \begin{array}{l}
    u_1 - u_1^+ = \lambda^k (a_{11}  C_k(\phi) + a_{12} \mu_2 S_k(\phi)) x_1 + \lambda^k (a_{11}  S_k(\phi) + a_{12} C_k(\phi)) x_2 + \\ 
    \qquad + b_1 (y - y^-)
    + \hat\lambda^{k} \tilde\xi_{1k}(x, y, \mu) +    
           O(\lambda^{2k}\|x\|^2 + \lambda^k \|x\| \cdot |y - y^-| + (y - y^-)^2), \\
    
    u_2 - u_2^+ = \lambda^k (a_{21}  C_k(\phi) + a_{22} \mu_2 S_k(\phi)) x_1 + \lambda^k (a_{21}  S_k(\phi) + a_{22} C_k(\phi)) x_2 + \\ 
    \qquad + b_2 (y - y^-) + \hat\lambda^{k} \tilde\xi_{2k}(x, y, \mu) +    
     O(\lambda^{2k}\|x\|^2 + \lambda^k \|x\| \cdot |y - y^-| + (y - y^-)^2), \\

    \gamma_2^{-j} v +  \hat\gamma_2^{-j} \eta_{3j}(u, v, \mu) = 
    \lambda^k (c_{1}  C_k(\phi) + c_{2} \mu_2 S_k(\phi)) x_1 + \lambda^k (c_{1}  S_k(\phi) + c_{2} C_k(\phi)) x_2 + \\ 
    \qquad + d_1 (y - y^-) + \hat\lambda^{k} \tilde\xi_{3k}(x, y, \mu) +    
     O(\lambda^{2k}\|x\|^2 +\lambda^k \|x\| \cdot |y - y^-| + (y - y^-)^2),
     \end{array}
\end{equation}
\begin{equation}\label{eq:Tkj121a}     
     \begin{array}{l}
    \bar x_1 - x_1^+ = a_{31}\nu_1^j u_1 + \hat\nu^{j} \tilde\eta_{1j}(u, v, \mu) +
        b_3 (v - v^-) + 
        O(\nu_1^{2j} \|u\|^2 + \nu_1^j |u| \cdot |v - v^-| + (v - v^-)^2), \\

    \bar x_2 - x_2^+ = a_{41}\nu_1^j u_1 + \hat\nu^{j} \tilde\eta_{2j}(u, v, \mu) +
        b_4 (v - v^-) + 
        O(\nu_1^{2j} \|u\|^2 + \nu_1^j |u| \cdot |v - v^-| + (v - v^-)^2), \\ 

    \gamma_1^{-k} \bar y +  \hat\gamma_1^{-k} \xi_{2k}(\bar x, \bar y, \mu) = \mu_1 + 
    c_3 \nu_1^j u_1 + \hat\nu^j \tilde\eta_{3j}(u, v, \mu) +     
         d_2 (v - v^-)^2 + \\ \qquad + O(\nu_1^{2j} \|u\|^2 + \nu_1^j \|u\| \cdot |v - v^-| +  (v - v^-)^3).
  \end{array}
\end{equation}

Make a coordinate shift
$u_{new} = u - u^+ + \varphi_{kj}^1$, $v_{new} = v - v^- +
\varphi_{kj}^2$, $x_{new} = x - x^+ + \psi_{kj}^1$, $y_{new} = y - y^- + \psi_{kj}^2$,
where $\varphi_{kj}^1,\; \psi_{kj}^2 = O(\gamma_2^{-j} + k\lambda^k)$ and 
$\varphi_{kj}^2,\;  \psi_{kj}^1 = O(\nu_1^j)$. With that, the nonlinearity functions in the left parts of the third equations 
of (\ref{eq:Tkj12a}) and (\ref{eq:Tkj121a}) can be expressed as Taylor expansions 
$\eta_{3j}(u + u^+ + \varphi_{kj}^1, v + v^- + \varphi_{kj}^2, \mu)) = 
\eta_{3j}^0 + \eta_{3j}^1 v + \eta_{3j}^2(u, v) + \eta_{3j}^3(v)$, $\xi_{2k}(\bar x + x^+ + \psi_{ij}^1, \bar y + y^- + \psi_{ij}^2, \mu)) = 
\xi_{2k}^0 + \xi_{2k}^1 \bar y + \xi_{2k}^2(\bar x, \bar y) + \xi_{2k}^3(\bar y)$ respectively, where coefficients  
$\eta_{3j}^0$, $\eta_{3j}^1$, $\xi_{2k}^0$, $\xi_{2k}^1$ are uniformly bounded in $k$ and $j$ for all small $\mu$ and
$\eta_{3j}^2 (u, v) = O(u)$, $\xi_{2k}^2 (\bar x, \bar y) = O(\bar x)$, $\eta_{3j}^3(v) = O(v^2)$, $\xi_{2k}^3(\bar y) = O(\bar y^2)$. 
We select constants $\varphi_{kj}^1, \; \varphi_{kj}^2,\; \psi_{kj}^1,\; \psi_{kj}^2$ in such a way
that all constant terms in equations (\ref{eq:Tkj12a}), the constant terms in the first two equations and the linear in $v_{new}$ term
in the last equation of (\ref{eq:Tkj121a}) vanish. In addition, we plug the expressions for $u$ coordinates from the
first two equations of (\ref{eq:Tkj12a}) into the third one, this will cause additions of order $\hat\gamma_2^{-j}$ to all
the coefficients. The system is rewritten as:
\begin{equation}\label{eq:3a}
    \begin{array}{l}
    u_1  = \lambda^k (a_{11}  C_k(\phi) + a_{12} \mu_2 S_k(\phi)) x_1 + \lambda^k (a_{11}  S_k(\phi) + a_{12} C_k(\phi)) x_2 + b_1 y  + \\ 
    \qquad 
           + O(\hat\lambda^{k}\|x\|^2 + \lambda^k \|x\| \cdot |y| + 
           y^2), \\
    
    u_2 = \lambda^k (a_{21}  C_k(\phi) + a_{22} \mu_2 S_k(\phi)) x_1 + \lambda^k (a_{21}  S_k(\phi) + a_{22} C_k(\phi)) x_2 + b_2 y  +\\ 
    \qquad +     
     O(\hat\lambda^{k}\|x\|^2 + \lambda^k \|x\| \cdot |y| + y^2), \\

    \gamma_2^{-j}(1 + q_{kj}^2) v +  \hat\gamma_2^{-j} O(v^2) = 
    \lambda^k (c_{1}  C_k(\phi) + c_{2} \mu_2 S_k(\phi)) x_1 + \lambda^k (c_{1}  S_k(\phi) + c_{2} C_k(\phi)) x_2 + \\ 
    \qquad + d_1 y +    
     O(\hat\lambda^{k}\|x\|^2 +\lambda^k \|x\| \cdot |y | + y^2),
     \end{array}
\end{equation}
\begin{equation}\label{eq:4a}     
     \begin{array}{l}
    \bar x_1 = a_{31} \nu_1^j u_1 + \tilde a_{32} \hat\nu^j u_2 + b_3 v +  O( \hat\nu^j \|u\|^2 + \nu_1^j \|u\| \cdot |v| + v^2), \\

    \bar x_2 = a_{41} \nu_1^j u_1 + \tilde a_{42} \hat\nu^j u_2 + b_4 v + O(\hat\nu^j \|u\|^2 + \nu_1^j \|u\| \cdot |v| + v^2), \\

    \gamma_1^{-k} (1 + q_{kj}^{1})\bar y + \hat\gamma_1^{-k} O(\bar x) + \hat\gamma_1^{-k} O(\bar y^2) = M^1 + c_3 \nu_1^j u_1 + 
    \tilde c_4 \hat\nu^j u_2 + d_2 v^2 + \\ \qquad + O(\hat\nu^{j} \|u\|^2 + \nu_1^j \|u\| \cdot |v| + |v|^3),
  \end{array}
\end{equation}
where $q_{kj}^{1} = O\left(\left( \hat\gamma_1 / \gamma_1 \right)^{-k} \right)$, 
$q_{kj}^{2} = O\left(\left( \hat\gamma_2 / \gamma_2 \right)^{-j} \right)$, coefficients marked with ``tilde'' are uniformly bounded
for small $\mu$ and the following expression is valid for $M^1$:

\begin{equation} \label{M1a}
M^1 = \mu_1 + \nu_1^j c_3 u_1^+ - \gamma_1^{-k}y^- + O(\hat\gamma_1^{-k} + \hat\nu^j + \gamma_1^{-k}\gamma_2^{-j}). 
\end{equation}

Next, we take the right-hand side of the third equation of (\ref{eq:3a}) divided by the factor 
$(1 + q_{kj}^{2})$ from the left-hand side as the new variable $y$ --  the equation becomes 
$\gamma_2^{-j} v +  \hat\gamma_2^{-j}  O(v^2) = y$, and we substitute it into the third equation of (\ref{eq:4a}). We also define $u_{1new} = u_1 - \frac{b_{1}}{d_1} \gamma_2^{-j} v + \hat\gamma_2^{-j} O(v^2)$ and $u_{2new} = u_2 - \frac{b_{2}}{d_1} \gamma_2^{-j} v + \hat\gamma_2^{-j} O(v^2)$ such that we remove all terms in the equations for $u$ which depend
on $y$ alone. 
In addition, we substitute the expressions for $\bar x$ to the last equation of (\ref{eq:4a}). These actions cause
the linear in $v$ term of order $O(\hat\gamma_1^{-k} + \nu_1^j \gamma_2^{-j})$ to appear in the equation for $\bar v$. We will
make it zero again later.
Thus, we obtain
\begin{equation}\label{eq:Tijshc}
    \begin{array}{l}

u_1 = \lambda^k (A_{11}  C_k(\phi) + A_{12} \mu_2 S_k(\phi)) x_1 + \lambda^k (A_{11}  S_k(\phi) + A_{12} C_k(\phi)) x_2  + \\ 
    \qquad 
     + O(\hat\lambda^{k}\|x\|^2 + \lambda^k\gamma_2^{-j} \|x\| \cdot |v|), \\
 
u_2 = \lambda^k (A_{21}  C_k(\phi) + A_{22} \mu_2 S_k(\phi)) x_1 + \lambda^k (A_{21}  S_k(\phi) + A_{22} C_k(\phi)) x_2  + \\ 
    \qquad 
     + O(\hat\lambda^{k}\|x\|^2 + \lambda^k\gamma_2^{-j} \|x\| \cdot |v|), \\

\bar x_1 = a_{31} \nu_1^j u_1 + \tilde a_{32} \hat\nu^j u_2 + b_3 v +
    O(\hat\nu^j \|u\|^2 + \nu_1^j\|u\| \cdot |v| + v^2), \\

\bar x_2 = a_{41} \nu_1^j u_1 + \tilde a_{42} \hat\nu^j u_2 + b_4 v +
    O(\hat\nu^j \|u\|^2 + \nu_1^{j}\|u\| \cdot |v| + v^2), \\

\displaystyle  \frac{\gamma_1^{-k}\gamma_2^{-j}}{d_1} \bar v(1 + q_{kj}^{3}) + \gamma_1^{-k}\hat\gamma_2^{-j} O(\bar v^2) =  M^1 + c_3 \nu_1^j u_1 + 
    \tilde c_4 \hat\nu^j u_2 + \\ 
      \qquad\qquad\qquad\qquad  + O(\hat\gamma_1^{-k} + \nu_1^j \gamma_2^{-j}) v 
      + d_2 v^2 + O(\hat\nu^j \|u\|^2 + \nu_1^{j}\|u\| \cdot |v| + |v|^3),
  \end{array}
\end{equation}
where $q_{ij}^{3} = O\left(\left( \hat\gamma_1 / \gamma_1 \right)^{-k} + \left(\hat\gamma_2 / \gamma_2 \right)^{-j} \right)$ and 
\begin{equation} \label{eq:Anew}
\begin{array}{c}
A_{11} = a_{11} - b_1 c_1 / d_1, \; A_{12} = a_{12} - b_1 c_2 / d_1 \\
A_{21} = a_{21} - b_2 c_1 / d_1, \; A_{22} = a_{22} - b_2 c_2 / d_1.
\end{array}
\end{equation}

We substitute $u$ as a function of $x$ and $v$ from the first two equations into the last three ones. After this, in addition, 
we make a shift of $(x, v)$ coordinates by a constant of order $O(\hat\gamma_1^{-k} + \nu_1^j \gamma_2^{-j})$ to 
nullify the linear in $v$ term in the last equation. This gives
us the following formula in variables $(x, v)$:
\begin{equation}\label{eq:Tij4}
  \begin{array}{l}
    \bar x_1 =  a_{31}\lambda^k \nu_1^j (A_{11}  C_k(\phi) + A_{12} \mu_2 S_k(\phi)) x_1 + a_{31} \lambda^k \nu_1^j (A_{11}  S_k(\phi) + A_{12} C_k(\phi)) x_2 + 
    b_3 v + \\ 
    \qquad\qquad + O(\hat\lambda^k \nu_1^j + \lambda^k \nu_1^j \gamma_2^{-j}) O(\|x\|^2) + 
    \lambda^k \nu_1^j O(\|x\| \cdot |v|) + O(v^2), \\

   \bar x_2 =  a_{41}\lambda^k \nu_1^j (A_{11}  C_k(\phi) + A_{12} \mu_2 S_k(\phi)) x_1 + a_{41} \lambda^k \nu_1^j (A_{11}  S_k(\phi) + A_{12} C_k(\phi)) x_2 + 
    b_4 v + \\ 
    \qquad\qquad + O(\hat\lambda^k \nu_1^j + \lambda^k \nu_1^j \gamma_2^{-j}) O(\|x\|^2) + 
    \lambda^k \nu_1^j O(\|x\| \cdot |v|) + O(v^2), \\ 

    \displaystyle  \frac{\gamma_1^{-k}\gamma_2^{-j}}{d_1} \bar v(1 + q_{kj}^{3}) + \gamma_1^{-k}\hat\gamma_2^{-j} O(\bar v^2) = 
    M^1 +  c_3\lambda^k \nu_1^j (A_{11}  C_k(\phi) + A_{12} \mu_2 S_k(\phi)) x_1 + \\ 
    \qquad\qquad + c_3 \lambda^k \nu_1^j (A_{11}  S_k(\phi) + A_{12} C_k(\phi)) x_2  + d_2 v^2 + \\

    \qquad\qquad + O(\hat\lambda^k \nu_1^j + \lambda^k \nu_1^j \gamma_2^{-j}) O(\|x\|^2) + \lambda^k \nu_1^j O(\|x\| \cdot |v|) + O(v^3),
 \end{array}
\end{equation}

In order to eliminate the terms that depend on $v$ only in the first equation, we make the following coordinate change: $x_{1new} = x_1 - b_3/b_4 x_2$. The first return map takes the form:
\begin{equation}\label{TkjIIa}
    \begin{array}{l}
         \bar x_1 = \tilde A_{11}(\phi)\lambda^k \nu_1^j x_1 + \tilde A_{12}(\phi)\lambda^k \nu_1^j x_2 +   O(\hat\lambda^k \nu_1^j + \lambda^k \nu_1^j \gamma_2^{-j}) O(\|x\|^2) + 
    \lambda^k \nu_1^j O(\|x\| \cdot |v|) \\

    \bar x_2 = \tilde A_{21}(\phi)\lambda^k \nu_1^j x_1 + \tilde A_{22}(\phi)\lambda^k \nu_1^j x_2 + b_4 v + O(\hat\lambda^k \nu_1^j + \lambda^k \nu_1^j \gamma_2^{-j}) O(\|x\|^2) + 
    \lambda^k \nu_1^j O(\|x\| \cdot |v|) + O(v^2) \\

     \displaystyle  \frac{\gamma_1^{-k}\gamma_2^{-j}}{d_1} \bar v(1 + q_{kj}^{3}) + \gamma_1^{-k}\hat\gamma_2^{-j} O(\bar v^2) = 
    M^1 +  \tilde C_{1}(\phi)\lambda^k \nu_1^j x_1 
    + \tilde C_2 \lambda^k \nu_1^j x_2  + d_2 v^2 + \\

    \qquad\qquad + O(\hat\lambda^k \nu_1^j + \lambda^k \nu_1^j \gamma_2^{-j}) O(\|x\|^2) + \lambda^k \nu_1^j O(\|x\| \cdot |v|) + O(v^3),
          
    \end{array}
\end{equation}
where
\begin{equation}
\begin{array}{l}
\displaystyle \tilde A_{11}(\phi) = 
\left( a_{31} - \frac{b_3}{b_4} a_{41}\right) (A_{11}C_k(\phi) + A_{12}\mu_2 S_k(\phi)) + O\left(\left(\hat\lambda/\lambda\right)^k  +  \left(\hat\nu/\nu_1\right)^j\right)
\\
\displaystyle \tilde A_{12}(\phi) = 
\frac{b_3}{b_4} \left( a_{31} - \frac{b_3}{b_4} a_{41}\right) (A_{11}C_k(\phi) + A_{12}\mu_2 S_k(\phi)) + \\
\displaystyle \qquad + \left( a_{31} - \frac{b_3}{b_4} a_{41}\right) (A_{11}S_k(\phi) + A_{12}C_k(\phi)) +
O\left(\left(\hat\lambda/\lambda\right)^k  +  \left(\hat\nu/\nu_1\right)^j\right)
\\
       
\displaystyle \tilde A_{21}(\phi) = a_{41} (A_{11}C_k(\phi) + A_{12}\mu_2 S_k(\phi)) + O\left(\left(\hat\lambda/\lambda\right)^k  +  \left(\hat\nu/\nu_1\right)^j\right)\\
\displaystyle \tilde A_{22}(\phi) = 
a_{41} \left( A_{11}C_k(\phi) + A_{12}\mu_2 S_k(\phi) + 
\frac{b_3}{b_4}\left( A_{11}S_k(\phi) + A_{12}C_k(\phi)\right)\right) +\\
\qquad + O\left(\left(\hat\lambda/\lambda\right)^k  +  \left(\hat\nu/\nu_1\right)^j\right) \\
\displaystyle \tilde C_{1}(\phi) = c_3 (A_{11}C_k(\phi) + A_{12}\mu_2 S_k(\phi))  + O\left(\left(\hat\lambda/\lambda\right)^k  +  \left(\hat\nu/\nu_1\right)^j\right) \\
\displaystyle \tilde  C_{2}(\phi) = c_3 \left( A_{11}C_k(\phi) + A_{12}\mu_2 S_k(\phi) + 
\frac{b_3}{b_4}\left( A_{11}S_k(\phi) + A_{12}C_k(\phi)\right)\right) +\\
\qquad + O\left(\left(\hat\lambda/\lambda\right)^k  +  \left(\hat\nu/\nu_1\right)^j\right).
\end{array}
\end{equation}

Take $\mu_2 < 0$. Coefficient $\tilde C_2(\phi)$ can be then represented as 
$\tilde C_2(\phi) = L_1(\mu) \cos k \phi + L_2(\mu) \sin k \phi = L(\mu) \sin (k \phi + \psi)$. We will vary $\mu_2$ near those values, for which $\tilde C_{2}(\phi)$ vanishes, so that it becomes asymptotically small for $k,j \to \infty$. Note that $\tilde C_1(\phi)$ and $\tilde A_{12}(\phi)$ will remain bounded from zero at these parameter values, because otherwise the first return map will be not a diffeomorphism. Finally, we rescale the coordinates as follows:
\begin{equation}
\begin{array}{c}
\displaystyle  v = -\frac{\gamma_1^{-k}\gamma_2^{-j}}{d_{1} d_{2}} (1 + q_{kj}^{3}) \; Y \;,\; 
x_1  = - \frac{b_4 \tilde A_{12}  \lambda^k \nu_1^j\gamma_1^{-k} \gamma_2^{-j}} {d_{1} d_{2}}(1 + q_{kj}^{3})\;X_1 \;,\; 
x_2  = - \frac{b_4 \gamma_1^{-k}  \gamma_2^{-j}}{d_{1} d_{2}}(1 + q_{kj}^{3})\;X_2 \;.
\end{array}
\end{equation}.

The resulting map is asymptotically close to three-dimensional Henon map (\ref{H3D}) as $k,j \to \infty$, and formulas (\ref{M123IIa}) are valid for the rescaled parameters. For small variations of $\mu_1$ parameter $M_1$ takes arbitrary finite values. We vary $\mu_2 < 0$ near a zero of function 
$\sin (k \phi + \psi)$, and being multiplied by $(\lambda \gamma_1)^k (\nu_1 \gamma_2)^j \sim \lambda^{-k}\nu_2^{-j}$ we obtain that parameter $M_2$ also takes arbitrary finite values. We also select arbitrary large $k$ and $j$ in such a way that the Jacobian of the first return map $B = J_{12} J_{21} (\lambda^2 \gamma_1)^k (\nu_1 \nu_2 \gamma_2)^j$ stays finite. Then, by small variations of $\mu_3$ it will also take arbitrary finite values, positive or negative.

\subsection{Heteroclinic Case IIb.}\label{sect:caseIIb}

Now we construct the first return map substituting local maps  (\ref{eq:T0kk1}) and (\ref{eq:T20j}) to  global maps 
(\ref{eq:T12b}) and 
(\ref{eq:T21b}):
\begin{equation}\label{eq:Tkj12b}
    \begin{array}{l}
    u_1 - u_1^+ = \lambda^k (a_{11}  C_k(\phi) + a_{12} \mu_2 S_k(\phi)) x_1 + \lambda^k (a_{11}  S_k(\phi) + a_{12} C_k(\phi)) x_2 + \\ 
    \qquad + b_1 (y - y^-)
    + \hat\lambda^{k} \tilde\xi_{1k}(x, y, \mu) +    
           O(\lambda^{2k}\|x\|^2 + \lambda^k \|x\| \cdot |y - y^-| + (y - y^-)^2), \\
    
    u_2 - u_2^+ = \lambda^k (a_{21}  C_k(\phi) + a_{22} \mu_2 S_k(\phi)) x_1 + \lambda^k (a_{21}  S_k(\phi) + a_{22} C_k(\phi)) x_2 + \\ 
    \qquad + b_2 (y - y^-) + \hat\lambda^{k} \tilde\xi_{2k}(x, y, \mu) +    
     O(\lambda^{2k}\|x\|^2 + \lambda^k \|x\| \cdot |y - y^-| + (y - y^-)^2), \\

    \gamma_2^{-j} v +  \hat\gamma_2^{-j} \eta_{3j}(u, v, \mu) = \mu_1 +
    \lambda^k (c_{1}  C_k(\phi) + c_{2} \mu_2 S_k(\phi)) x_1 + \lambda^k (c_{1}  S_k(\phi) + c_{2} C_k(\phi)) x_2 + \\ 
    \qquad + d_1 (y - y^-)^2 + \hat\lambda^{k} \tilde\xi_{3k}(x, y, \mu) +    
     O(\lambda^{2k}\|x\|^2 +\lambda^k \|x\| \cdot |y - y^-| + (y - y^-)^3),
     \end{array}
\end{equation}
\begin{equation}\label{eq:Tkj121b}     
     \begin{array}{l}
    \bar x_1 - x_1^+ = a_{31}\nu_1^j u_1 + \hat\nu^{j} \tilde\eta_{1j}(u, v, \mu) +
        b_3 (v - v^-) + 
        O(\nu_1^{2j} \|u\|^2 + \nu_1^j |u| \cdot |v - v^-| + (v - v^-)^2), \\

    \bar x_2 - x_2^+ = a_{41}\nu_1^j u_1 + \hat\nu^{j} \tilde\eta_{2j}(u, v, \mu) +
        b_4 (v - v^-) + 
        O(\nu_1^{2j} \|u\|^2 + \nu_1^j |u| \cdot |v - v^-| + (v - v^-)^2), \\ 

    \gamma_1^{-k} \bar y +  \hat\gamma_1^{-k} \xi_{2k}(\bar x, \bar y, \mu) = 
    c_3 \nu_1^j u_1 + \hat\nu^j \tilde\eta_{3j}(u, v, \mu) +     
         d_2 (v - v^-) + \\ \qquad + O(\nu_1^{2j} \|u\|^2 + \nu_1^j \|u\| \cdot |v - v^-| +  (v - v^-)^2).
  \end{array}
\end{equation}

Make a coordinate shift
$u_{new} = u - u^+ + \varphi_{kj}^1$, $v_{new} = v - v^- +
\varphi_{kj}^2$, $x_{new} = x - x^+ + \psi_{kj}^1$, $y_{new} = y - y^- + \psi_{kj}^2$,
where $\varphi_{kj}^1,\; \psi_{kj}^2 = O(\gamma_2^{-j} + k\lambda^k)$ and 
$\varphi_{kj}^2,\;  \psi_{kj}^1 = O(\nu_1^j)$. With that, the nonlinearity functions in the left parts of the third equations 
of (\ref{eq:Tkj12a}) and (\ref{eq:Tkj121a}) can be expressed as Taylor expansions 
$\eta_{3j}(u + u^+ + \varphi_{kj}^1, v + v^- + \varphi_{kj}^2, \mu)) = 
\eta_{3j}^0 + \eta_{3j}^1 v + \eta_{3j}^2(u, v) + \eta_{3j}^3(v)$, $\xi_{2k}(\bar x + x^+ + \psi_{ij}^1, \bar y + y^- + \psi_{ij}^2, \mu)) = 
\xi_{2k}^0 + \xi_{2k}^1 \bar y + \xi_{2k}^2(\bar x, \bar y) + \xi_{2k}^3(\bar y)$ respectively, where coefficients  
$\eta_{3j}^0$, $\eta_{3j}^1$, $\xi_{2k}^0$, $\xi_{2k}^1$ are uniformly bounded in $k$ and $j$ for all small $\mu$ and
$\eta_{3j}^2 (u, v) = O(u)$, $\xi_{2k}^2 (\bar x, \bar y) = O(\bar x)$, $\eta_{3j}^3(v) = O(v^2)$, $\xi_{2k}^3(\bar y) = O(\bar y^2)$. 
We select constants $\varphi_{kj}^1, \; \varphi_{kj}^2,\; \psi_{kj}^1,\; \psi_{kj}^2$ in such a way
that all constant terms in equations (\ref{eq:Tkj121b}), the constant terms in the first two equations and the linear in $y_{new}$ term
in the last equation of (\ref{eq:Tkj12b}) vanish. In addition, we plug the expressions for $\bar x$ coordinates from the
first two equations of (\ref{eq:Tkj121b}) into the third one, this will cause additions of order $\hat\gamma_1^{-k}$ to all
the coefficients. The system is rewritten as:
\begin{equation}\label{eq:3b}
    \begin{array}{l}
    u_1  = \lambda^k (a_{11}  C_k(\phi) + a_{12} \mu_2 S_k(\phi)) x_1 + \lambda^k (a_{11}  S_k(\phi) + a_{12} C_k(\phi)) x_2 + b_1 y  + \\ 
    \qquad 
           + O(\hat\lambda^{k}\|x\|^2 + \lambda^k \|x\| \cdot |y| + 
           y^2), \\
    
    u_2 = \lambda^k (a_{21}  C_k(\phi) + a_{22} \mu_2 S_k(\phi)) x_1 + \lambda^k (a_{21}  S_k(\phi) + a_{22} C_k(\phi)) x_2 + b_2 y  +\\ 
    \qquad +     
     O(\hat\lambda^{k}\|x\|^2 + \lambda^k \|x\| \cdot |y| + y^2), \\

    \gamma_2^{-j}(1 + q_{kj}^2) v + \hat\gamma_2^{-j} O(u) 
 + \hat\gamma_2^{-j} O(v^2) = M^1  +
    \lambda^k (c_{1}  C_k(\phi) + c_{2} \mu_2 S_k(\phi)) x_1 + \\ 
    \qquad + \lambda^k (c_{1}  S_k(\phi) + c_{2} C_k(\phi)) x_2  
     + d_1 y^2 +    
     O(\hat\lambda^{k}\|x\|^2 +\lambda^k \|x\| \cdot |y | + |y|^3),
     \end{array}
\end{equation}
\begin{equation}\label{eq:4b}     
     \begin{array}{l}
    \bar x_1 = a_{31} \nu_1^j u_1 + \tilde a_{32} \hat\nu^j u_2 + b_3 v +  O( \hat\nu^j \|u\|^2 + \nu_1^j \|u\| \cdot |v| + v^2), \\

    \bar x_2 = a_{41} \nu_1^j u_1 + \tilde a_{42} \hat\nu^j u_2 + b_4 v + O(\hat\nu^j \|u\|^2 + \nu_1^j \|u\| \cdot |v| + v^2), \\

    \gamma_1^{-k} (1 + q_{kj}^{1})\bar y  + \hat\gamma_1^{-k} O(\bar y^2) =  c_3 \nu_1^j u_1 + 
    \tilde c_4 \hat\nu^j u_2 + d_2 v + \\ \qquad + O(\hat\nu^{j} \|u\|^2 + \nu_1^j \|u\| \cdot |v| + v^2),
  \end{array}
\end{equation}
where $q_{kj}^{1} = O\left(\left( \hat\gamma_1 / \gamma_1 \right)^{-k} \right)$, 
$q_{kj}^{2} = O\left(\left( \hat\gamma_2 / \gamma_2 \right)^{-j} \right)$, coefficients marked with ``tilde'' are uniformly bounded
for small $\mu$ and the following expression is valid for $M^1$:

\begin{equation} \label{M1b}
M^1 = \mu_1 + \lambda^k( (c_1 C_k(\phi) + c_2 \mu_2 S_k(\phi))x_1^+  + 
(c_1 S_k(\phi) + c_2 C_k(\phi))x_2^+) - \gamma_2^{-j}v^- + O(\hat\gamma_2^{-j} + \hat\lambda^k + \gamma_1^{-k}\gamma_2^{-j}). 
\end{equation}

Next, we take the right-hand side of the third equation of (\ref{eq:4b}) divided by the factor 
$(1 + q_{kj}^{1})$ from the left-hand side as the new variable $v$ --  the equation becomes 
$\gamma_1^{-k} \bar y + \hat\gamma_1^{-k}  O(\bar y^2) = v$, and we substitute it into the third equation of (\ref{eq:3b}). We also define $x_{1new} = x_1 - \frac{b_{3}}{d_2} \gamma_1^{-k} y + \hat\gamma_1^{-k} O(y^2)$ and $x_{2new} = x_2 - \frac{b_{4}}{d_2} \gamma_1^{-k} y + \hat\gamma_1^{-k} O(y^2)$ such that we remove all terms in the equations for $x$ which depend
on $v$ alone. 
In addition, we substitute the expressions for $\bar x$ to the last equation of (\ref{eq:4b}). These actions cause
the linear in $y$ term of order $O(\hat\gamma_2^{-j} + \lambda^k \gamma_1^{-k})$ to appear in the equation for $\bar y$. We will
make it zero again later.
Thus, we obtain
\begin{equation}\label{eq:Tijshcb}
    \begin{array}{l}

u_1 = \lambda^k (a_{11}  C_k(\phi) + a_{12} \mu_2 S_k(\phi)) x_1 + \lambda^k (a_{11}  S_k(\phi) + a_{12} C_k(\phi)) x_2  + b_1 v\\ 
    \qquad 
     + O(\hat\lambda^{k}\|x\|^2 + \lambda^k\gamma_2^{-j} \|x\| \cdot |v| + v^2), \\
 
u_2 = \lambda^k (a_{21}  C_k(\phi) + a_{22} \mu_2 S_k(\phi)) x_1 + \lambda^k (a_{21}  S_k(\phi) + a_{22} C_k(\phi)) x_2  + b_2 v\\ 
    \qquad 
     + O(\hat\lambda^{k}\|x\|^2 + \lambda^k\gamma_2^{-j} \|x\| \cdot |v| + v^2), \\

\bar x_1 = A_{31} \nu_1^j u_1 + \tilde A_{32} \hat\nu^j u_2  +
    O(\hat\nu^j \|u\|^2 + \nu_1^j\|u\| \cdot |v| ), \\

\bar x_2 = A_{41} \nu_1^j u_1 + \tilde A_{42} \hat\nu^j u_2 + 
    O(\hat\nu^j \|u\|^2 + \nu_1^{j}\|u\| \cdot |v| ), \\

\displaystyle  \frac{\gamma_1^{-k}\gamma_2^{-j}}{d_2} \bar y(1 + q_{kj}^{3}) + \hat\gamma_1^{-k} \gamma_2^{-j} O(\bar y^2) =  M^1 + \lambda^k (c_{1}  C_k(\phi) + c_{2} \mu_2 S_k(\phi)) x_1 + \\ 
    \qquad + \lambda^k (c_{1}  S_k(\phi) + c_{2} C_k(\phi)) x_2 +  
      O(\hat\gamma_2^{-j} + \lambda^k \gamma_1^{-k}) y 
      + d_1 y^2 + \\
      \qquad + O(\hat\lambda^k \|x\|^2 + \lambda^{k}\|x\| \cdot |y| + |y|^3),
  \end{array}
\end{equation}
where $q_{ij}^{3} = O\left(\left( \hat\gamma_1 / \gamma_1 \right)^{-k} + \left(\hat\gamma_2 / \gamma_2 \right)^{-j} \right)$ and 
%
$
A_{31} = a_{31} - b_3 c_3 / d_2, \; 
A_{41} = a_{41} - b_4 c_3 / d_2. \; 
$

We substitute $u$ as a function of $x$ and $y$ from the first two equations into the last three ones. After this, in addition, 
we make a shift of $(x, y)$ coordinates by a constant of order $O(\hat\gamma_2^{-j} + \lambda^k \gamma_1^{-k})$ to 
nullify the linear in $y$ term in the last equation. This gives
us the following formula in variables $(x, y)$:
\begin{equation}\label{eq:Tij4b}
  \begin{array}{l}
    \bar x_1 =  A_{31}\lambda^k \nu_1^j (a_{11}  C_k(\phi) + a_{12} \mu_2 S_k(\phi)) x_1 + A_{31} \lambda^k \nu_1^j (a_{11}  S_k(\phi) + a_{12} C_k(\phi)) x_2 + 
    A_{31} b_1\nu_1^j y + \\ 
    \qquad\qquad + O( \lambda^k \hat\nu^j + \lambda^k \nu_1^j \gamma_1^{-k}) O(\|x\|^2) + 
    \lambda^k \nu_1^j O(\|x\| \cdot |y|) + O(y^2), \\

   \bar x_2 =  A_{41}\lambda^k \nu_1^j (a_{11}  C_k(\phi) + a_{12} \mu_2 S_k(\phi)) x_1 + A_{41} \lambda^k \nu_1^j (a_{11}  S_k(\phi) + a_{12} C_k(\phi)) x_2 + 
    A_{41}b_1 \nu_1^j y + \\ 
    \qquad\qquad + O(\lambda^k \hat\nu^j + \lambda^k \nu_1^j \gamma_1^{-k}) O(\|x\|^2) + 
    \lambda^k \nu_1^j O(\|x\| \cdot |y|) + O(y^2), \\ 

    \displaystyle  \frac{\gamma_1^{-k}\gamma_2^{-j}}{d_2} \bar y(1 + q_{kj}^{3}) + \hat\gamma_1^{-k}\gamma_2^{-j} O(\bar y^2) = 
    M^1 + \lambda^k (c_{1}  C_k(\phi) + c_{2} \mu_2 S_k(\phi)) x_1 + \\ 
    \qquad + \lambda^k (c_{1}  S_k(\phi) + c_{2} C_k(\phi)) x_2  + d_1 y^2 + \\

    \qquad\qquad + O(\lambda^k \hat\nu^j + \lambda^k \nu_1^j \gamma_1^{-k}) O(\|x\|^2) + \lambda^k \nu_1^j O(\|x\| \cdot |y|) + O(y^3),
 \end{array}
\end{equation}

In order to eliminate the terms that depend on $y$ only in the first equation, we make the following coordinate change: $x_{1new} = x_1 - A_{31}/A_{41} x_2$. The first return map takes the form:
\begin{equation}\label{TkjIIb}
    \begin{array}{l}
         \bar x_1 = \tilde A_{11}(\phi)\lambda^k \nu_1^j x_1 + \tilde A_{12}(\phi)\lambda^k \nu_1^j x_2 + 
  O( \lambda^k \hat\nu^j + \lambda^k \nu_1^j \gamma_1^{-k}) O(\|x\|^2) + 
    \lambda^k \nu_1^j O(\|x\| \cdot |y|)  \\

    \bar x_2 = \tilde A_{21}(\phi)\lambda^k \nu_1^j x_1 + \tilde A_{22}(\phi)\lambda^k \nu_1^j x_2 +       A_{41} b_1\nu_1^j y + O( \lambda^k \hat\nu^j + \lambda^k \nu_1^j \gamma_1^{-k}) O(\|x\|^2) + 
    \lambda^k \nu_1^j O(\|x\| \cdot |y|) + O(y^2) \\

     \displaystyle  \frac{\gamma_1^{-k}\gamma_2^{-j}}{d_2} \bar y(1 + q_{kj}^{3}) + \gamma_1^{-k}\hat\gamma_2^{-j} O(\bar y^2) = 
    M^1 +  \tilde C_{1}(\phi)\lambda^k  x_1 
    + \tilde C_2 \lambda^k  x_2  + d_1 y^2 + \\

    \qquad\qquad +  O( \lambda^k \hat\nu^j + \lambda^k \nu_1^j \gamma_1^{-k}) O(\|x\|^2) + 
    \lambda^k \nu_1^j O(\|x\| \cdot |y|) + O(y^3)
          
    \end{array}
\end{equation}
where
\begin{equation}
\begin{array}{l}
\displaystyle \tilde A_{11}(\phi) = 
O\left(\left(\hat\lambda/\lambda\right)^k  +  \left(\hat\nu/\nu_1\right)^j\right)
\\
\displaystyle \tilde A_{12}(\phi) = 
O\left(\left(\hat\lambda/\lambda\right)^k  +  \left(\hat\nu/\nu_1\right)^j\right)
\\
       
\displaystyle \tilde A_{21}(\phi) = 
A_{41}(a_{11}  C_k(\phi) + a_{12} \mu_2 S_k(\phi)) + O\left(\left(\hat\lambda/\lambda\right)^k  +  \left(\hat\nu/\nu_1\right)^j\right)\\
\displaystyle \tilde A_{22}(\phi) = 
A_{31}(a_{11}  C_k(\phi) + a_{12} \mu_2 S_k(\phi)) + A_{41} (a_{11}  S_k(\phi) + a_{12} C_k(\phi)) +
O\left(\left(\hat\lambda/\lambda\right)^k  +  \left(\hat\nu/\nu_1\right)^j\right) \\
\displaystyle \tilde C_{1}(\phi) = c_{1}  C_k(\phi) + c_{2} \mu_2 S_k(\phi) +O\left(\left(\hat\lambda/\lambda\right)^k  +  \left(\hat\nu/\nu_1\right)^j\right) \\
\displaystyle \tilde  C_{2}(\phi) = c_{1}  S_k(\phi) + c_{2} C_k(\phi) + \frac{A_{31}}{A_{41}}(c_{1}  C_k(\phi) + c_{2} \mu_2 S_k(\phi))   + O\left(\left(\hat\lambda/\lambda\right)^k  +  \left(\hat\nu/\nu_1\right)^j\right).
\end{array}
\end{equation}

Take $\mu_2 < 0$. Coefficient $\tilde C_2(\phi)$ can be then represented as 
$\tilde C_2(\phi) = L_1(\mu) \cos k \phi + L_2(\mu) \sin k \phi = L(\mu) \sin (k \phi + \psi)$. We will vary $\mu_2$ near those values, for which $\tilde C_{2}(\phi)$ vanishes so that it becomes asymptotically small for $k,j \to \infty$. Note that $\tilde C_1(\phi)$ and $\tilde A_{12}(\phi)$ will remain bounded from zero at these parameter values because otherwise, the first return map will be not a diffeomorphism. Finally, we rescale the coordinates as follows:
\begin{equation}
\begin{array}{c}
\displaystyle  y = -\frac{\gamma_1^{-k}\gamma_2^{-j}}{d_{1} d_{2}} (1 + q_{kj}^{3}) \; Y \;,\; 
x_1  = - \frac{b_1 A_{41} \tilde A_{12} \lambda^k \nu_1^{2j}\gamma_1^{-k} \gamma_2^{-j}} {d_{1} d_{2}}(1 + q_{kj}^{3})\;X_1 \\ 
\displaystyle x_2  = - \frac{b_1 A_{41} \nu_1^j \gamma_1^{-k}  \gamma_2^{-j}}{d_{1} d_{2}}(1 + q_{kj}^{3})\;X_2 \;.
\end{array}
\end{equation}.

The resulting map is asymptotically close to three-dimensional Henon map (\ref{H3D}) as $k,j \to \infty$, and formulas (\ref{M123IIa}) are valid for the rescaled parameters. For small variations of $\mu_1$ parameter $M_1$ takes arbitrary finite values. We vary $\mu_2 < 0$ near a zero of function 
$\sin (k \phi + \psi)$, and being multiplied by $(\lambda \gamma_1)^k (\nu_1 \gamma_2)^j \sim \lambda^{-k}\nu_2^{-j}$ we obtain that parameter $M_2$ also takes arbitrary finite values. We also select arbitrary large $k$ and $j$ in such a way that the Jacobian of the first return map $B = J_{12} J_{21} (\lambda^2 \gamma_1)^k (\nu_1 \nu_2 \gamma_2)^j$ stays finite. Then, by small variations of $\mu_3$ it will also take arbitrary finite values, positive or negative.


\section*{Acknowledgements}
The author thanks Sergey Gonchenko for the formulation of the problem and useful discussions, including the normal forms near saddles having eigenvalues with multiplicity greater than one.

\section*{Funding}
This paper is a contribution to the project M7 (Dynamics of Geophysical Problems in Turbulent Regimes) of the Collaborative Research Centre TRR 181 ``Energy Transfer in Atmosphere and Ocean'' funded by the Deutsche Forschungsgemeinschaft (DFG, German Research Foundation) -- Projektnummer 274762653. The paper is also supported by the
grant of the Russian Science
Foundation 19-11-00280.

\section*{Data availability statement}
The original contributions presented in the study are included
in the article/supplementary material, further inquiries can be
directed to the corresponding author.

\end{document}